# APPROXIMATION IN HIGHER-ORDER SOBOLEV SPACES AND HODGE SYSTEMS

PIERRE BOUSQUET, EMMANUEL RUSS, YI WANG, PO-LAM YUNG

ABSTRACT. Let $d \geq 2$ be an integer, $1 \leq l \leq d-1$ and $\varphi$ be a differential $l$-form on $\mathbb{R}^d$ with $\dot{W}^{1,d}$ coefficients. It was proved by Bourgain and Brezis ([5, Theorem 5]) that there exists a differential $l$-form $\psi$ on $\mathbb{R}^d$ with coefficients in $L^\infty \cap \dot{W}^{1,d}$ such that $d\varphi = d\psi$. In the same work, Bourgain and Brezis also left as an open problem the extension of this result to the case of differential forms with coefficients in the higher order space $\dot{W}^{2,d/2}$ or more generally in the fractional Sobolev spaces $\dot{W}^{s,p}$ with $sp = d$. We give a positive answer to this question, provided that $d - \kappa \leq l \leq d - 1$, where $\kappa$ is the largest positive integer such that $\kappa < \min(p, d)$. The proof relies on an approximation result (interesting in its own right) for functions in $\dot{W}^{s,p}$ by functions in $\dot{W}^{s,p} \cap L^\infty$, even though $\dot{W}^{s,p}$ does not embed into $L^\infty$ in this critical case. The proofs rely on some techniques due to Bourgain and Brezis but the context of higher order and/or fractional Sobolev spaces creates various difficulties and requires new ideas and methods.

**Keywords:** Hodge systems, approximation in limiting Sobolev spaces, higher-order Sobolev spaces, Triebel-Lizorkin spaces, Littlewood-Paley inequalities.

## Contents

1. INTRODUCTION

1.1. **The main result.** For $k \in \mathbb{N}$ and $1 < p < \infty$, let $\dot{W}^{k,p}(\mathbb{R}^d)$ be the homogeneous Sobolev space on $\mathbb{R}^d$, that is the completion of $C_c^\infty(\mathbb{R}^d)$ under the norm $\|\partial^k f\|_{L^p(\mathbb{R}^d)}$. It is well-known that while $\dot{W}^{k,p}(\mathbb{R}^d)$ embeds continuously into $L^{kp/(d-kp)}(\mathbb{R}^d)$ when $kp < d$, the embedding fails when $kp = d$. In a ground-breaking paper [5] (see also [4]), Bourgain and Brezis found a remedy for this failure when $k = 1$ and $p = d$. They showed that for any $f \in \dot{W}^{1,d}(\mathbb{R}^d)$ and any $\delta > 0$, there exists $F \in \dot{W}^{1,d} \cap L^\infty(\mathbb{R}^d)$ and a constant $C_\delta > 0$ independent of $f$, such that

$$\text{(1.1)} \qquad \sum_{i=1}^{d-1} \|\partial_i(f - F)\|_{L^d} \leq \delta \|f\|_{\dot{W}^{1,d}},$$

and

$$\text{(1.2)} \qquad \|F\|_{L^\infty} + \|F\|_{\dot{W}^{1,d}} \leq C_\delta \|f\|_{\dot{W}^{1,d}}.$$

This should be thought of as a theorem describing how a function in $\dot{W}^{1,d}(\mathbb{R}^d)$ can be approximated by a function in $L^\infty(\mathbb{R}^d)$, and the failure of the embedding of $\dot{W}^{1,d}(\mathbb{R}^d)$ into $L^\infty(\mathbb{R}^d)$ makes this result rather non-trivial. In the same paper, Bourgain and Brezis also derived many important consequences of this approximation theorem. Among them, they proved that if $l$ belongs to the set $[\![1, d-1]\!]$ (namely, the set of integers $j$ such that $1 \leq j \leq d-1$) and $\varphi$ is a differential $l$-form on $\mathbb{R}^d$ with $\dot{W}^{1,d}$ coefficients, then there exists a differential $l$-form $\psi$ on $\mathbb{R}^d$ with $\dot{W}^{1,d} \cap L^\infty$ coefficients such that

$$d\psi = d\varphi.$$

While the approximation result of Bourgain and Brezis in $\dot{W}^{1,d}$ was extended in a subelliptic context (see [33]), to our best knowledge, no extension to higher order Sobolev spaces was given until now, even for the space $\dot{W}^{2,d/2}(\mathbb{R}^d)$. In the present work, inspired by [5, Open Problem 2], we extend (1.1) and (1.2) to higher-order Sobolev spaces. To give an idea of our result, consider the following approximation problem for $\dot{W}^{2,d/2}(\mathbb{R}^d)$. Let $f \in \dot{W}^{2,d/2}(\mathbb{R}^d)$ and $\delta > 0$. We look for an approximation $F \in \dot{W}^{2,d/2}(\mathbb{R}^d) \cap L^\infty(\mathbb{R}^d)$ of $f$, in the sense that $\|\partial_i\partial_j(f - F)\|_{L^{d/2}}$ should be small for the largest possible set of indexes $i, j$. What the optimal statement should be is not obvious. For instance, can we approximate all the second order derivatives but one of them? Or, if we impose that $i \in [\![1, d-1]\!]$ (in order words, we avoid $i = d$), can we take any $j \in [\![1, d]\!]$? Should we introduce further restrictions? As a consequence of Theorem 1.1 below, we can give the following answer: denoting by $\kappa$ the largest positive integer satisfying $\kappa < d/2$, there exists $F \in \dot{W}^{2,d/2}(\mathbb{R}^d) \cap L^\infty(\mathbb{R}^d)$ such that, for all $i \in [\![1, \kappa]\!]$ and all $j \in [\![1, d]\!]$,

$$\text{(1.3)} \qquad \|\partial_i\partial_j(f - F)\|_{L^d} \leq \delta \|f\|_{\dot{W}^{2,d/2}}$$

along with

$$\|F\|_{L^\infty} + \|F\|_{\dot{W}^{2,d/2}} \leq C_\delta \|f\|_{\dot{W}^{2,d/2}}.$$

As we shall see, the restriction on the values of $i$ in (1.3) is sufficient to extend the results of [5] about differential forms to the case of $\dot{W}^{2,d/2}$ coefficients.

More generally, in this paper, we give an extension of the above results to a range of critical Triebel-Lizorkin spaces $\dot{F}_q^{\alpha,p}(\mathbb{R}^d)$ that barely fail to embed into $L^\infty$. In particular, our results cover the higher order Sobolev spaces $\dot{W}^{k,d/k}(\mathbb{R}^d)$ where $k$ is an integer with $1 < k < d$ (note that $\dot{W}^{k,d/k}(\mathbb{R}^d) = \dot{F}_2^{k,d/k}(\mathbb{R}^d)$), and the Sobolev spaces $\dot{W}^{\alpha,d/\alpha}(\mathbb{R}^d)$ of fractional order $\alpha \in (0,1)$, giving an answer to the Open Problem 2 in [5].

Our main result can be stated as follows:

**Theorem 1.1.** *Let $\alpha > 0$ and $p, q \in (1, \infty)$ such that $\alpha p = d$. Let $\kappa$ be the largest positive integer that satisfies $\kappa < \min\{p, d\}$. Then, for every $\delta > 0$, there exists a constant $C_\delta > 0$ such that, for every $f \in \dot{F}_q^{\alpha,p}(\mathbb{R}^d)$, there exists $F \in \dot{F}_q^{\alpha,p} \cap L^\infty(\mathbb{R}^d)$ such that*

$$\sum_{i=1}^{\kappa} \|\partial_i(f - F)\|_{\dot{F}_q^{\alpha-1,p}} \leq \delta \|f\|_{\dot{F}_q^{\alpha,p}},$$

*and*

$$\|F\|_{L^\infty} + \|F\|_{\dot{F}_q^{\alpha,p}} \leq C_\delta \|f\|_{\dot{F}_q^{\alpha,p}}.$$

From Theorem 1.1, we derive:

**Theorem 1.2.** *Let $\alpha > 0$, $p, q \in (1, \infty)$ such that $\alpha p = d$ and $l \in [\![d - \kappa, d - 1]\!]$, where $\kappa$ is the largest positive integer such that $\kappa < \min(p, d)$. Let $\varphi \in \dot{F}_q^{\alpha,p}(\Lambda^l \mathbb{R}^d)$. There exists $\psi \in \dot{F}_q^{\alpha,p}(\Lambda^l \mathbb{R}^d) \cap L^\infty(\Lambda^l \mathbb{R}^d)$ such that*

$$d\psi = d\varphi$$

*and*

$$\|\psi\|_{L^\infty(\Lambda^l \mathbb{R}^d)} + \|\psi\|_{\dot{F}_q^{\alpha,p}(\Lambda^l \mathbb{R}^d)} \lesssim \|d\varphi\|_{\dot{F}_q^{\alpha-1,p}(\Lambda^l \mathbb{R}^d)}.$$

By $\dot{F}_q^{\alpha,p}(\Lambda^l \mathbb{R}^d)$, we mean the space of differential $l$-forms on $\mathbb{R}^d$, the coefficients of which belong to $\dot{F}_q^{\alpha,p}(\mathbb{R}^d)$ (see Definition 2.1 below). The above statement extends the main result in [6] which was restricted to the conditions $\kappa = 1$ (which amounts to solving the equation div $X = f$ with $f \in \dot{F}_q^{\alpha,p}$), $\alpha > 1/2$ and $p \geq q \geq 2$; see also the earlier papers by Maz'ya [16] and also Mironescu [18] when $\kappa = 1$, and $p = q = 2$.

If we do not require the solution $\psi$ in Theorem 1.2 to be in $\dot{F}_q^{\alpha,p}$, then the theorem can be deduced from Proposition 2.1 of Van Schaftingen [30], which has a very elegant and simple proof. This elementary approach introduced in [27, 26] has been exploited in various settings, see in particular Lanzani and Stein [15] and Mitrea and Mitrea [19]. In the case of the equation div$X = f$ with $f \in L^d$, an algorithmic construction of a solution $X \in L^\infty$ was proposed in [24]. Up to our knowledge however, even for this equation and for the Sobolev space $\dot{W}^{1,d}$, there is no simple argument to prove the existence of a solution $\psi$ which is both in $L^\infty$ and in $\dot{W}^{1,d}$.

Many extensions, applications and recent developments of the original results established by Bourgain and Brezis in [5, 4] are presented in the excellent overview by Van Schaftingen [32]. We only quote some of them here:

(1) more general (higher order) operators than the exterior derivative have been studied in Van Schaftingen [28, 29, 31],

(2) similar problems have been considered when the space $\mathbb{R}^d$ is replaced by more general domains: half-spaces in Amrouche and Nguyen [1], smooth domains with specific boundary conditions in Brezis and Van Schaftingen [8], homogeneous groups in Chanillo and Van Schaftingen [9], Wang and Yung [33], symmetric spaces in Chanillo, Van Schaftingen and Yung [12, 11], and $CR$ manifolds in Yung [34].
(3) related Hardy inequalities were established by Maz'ya [17] (see also [7]),
(4) further applications of this theory can be found in Chanillo and Yung [13] and in Chanillo, Van Schaftingen and Yung [10].

1.2. **An overview of the strategy.** Let us first briefly recall the strategy of Bourgain and Brezis in their proof of the approximation theorem of $\dot{W}^{1,d}(\mathbb{R}^d)$ (that is, (1.1) and (1.2) above), before we turn to the difficulties we must face in proving Theorem 1.1. First they observe that for $f \in \dot{W}^{1,d}(\mathbb{R}^d)$, the Littlewood-Paley projections $\Delta_j f$ are uniformly bounded in $\mathbb{R}^d$ for all $j$, by Bernstein's inequality:
$$\|\Delta_j f\|_{L^\infty(\mathbb{R}^d)} \leq C\|\nabla f\|_{L^d(\mathbb{R}^d)}.$$
By normalizing $f$, one may thus assume that $\|\Delta_j f\|_{L^\infty(\mathbb{R}^d)} \leq 1$ for all $j \in \mathbb{Z}$. As a result, to approximate $f = \sum_{j \in \mathbb{Z}} \Delta_j f = \sum_{j \in \mathbb{Z}} \Delta_j f \cdot 1$ by a bounded function $F$, one is tempted to set
$$F(x) = \sum_{j \in \mathbb{Z}} \Delta_j f(x) \prod_{j' > j} (1 - |\Delta_{j'} f(x)|),$$
which would be automatically bounded by a partition of unity identity (see Lemma 3.2 below). Of course this cannot work, for this construction does not distinguish between the "good" directions $\partial_1, \ldots, \partial_{d-1}$ from the "bad" direction $\partial_d$ (whereas (1.1) distinguishes those). Thus Bourgain and Brezis introduce an auxiliary function $\omega_j(x)$, which controls $\Delta_j f(x)$ in the sense that
$$|\Delta_j f| \leq \omega_j \leq \|\Delta_j f\|_{L^\infty(\mathbb{R}^d)},$$
while satisfying good derivative bounds such as
$$|\partial_i \omega_j| \leq C 2^{j-\sigma} \omega_j \quad \text{for } i = 1, \ldots, d-1, \quad \text{and} \quad |\partial_d \omega_j| \leq C 2^j \omega_j;$$
here $\sigma > 0$ is a large parameter only depending on $\delta$. These $\omega_j$'s are constructed in [5] by using a sup-convolution (where one takes a supremum instead of an integral in the definition of a convolution), namely:
$$\omega_j(x) = \sup_{y \in \mathbb{R}^d} |\Delta_j f(y)| e^{-2^j |x_d - y_d| - 2^{j-\sigma}|x' - y'|},$$
where $x' - y' = (x_1 - y_1, \ldots, x_{d-1} - y_{d-1})$. With this in hand, one may be tempted to define the approximating function $F$ by setting
$$F(x) = \sum_{j \in \mathbb{Z}} \Delta_j f(x) \prod_{j' > j} (1 - \omega_j(x)),$$
which again would be automatically bounded by a partition of unity identity, and which has a better chance of obeying estimate (1.1). It turns out that this is still not sufficient; indeed, if $F$ were such defined, then
$$f - F = \sum_j \omega_j \mu_j$$

for some functions $\mu_j$ given by

$$\mu_j(x) = \sum_{j'<j} \Delta_{j'} f(x) \prod_{j'<j''<j} (1 - \omega_{j''}(x)).$$

These $\mu_j$'s are pointwisely bounded by 1 under our normalization of $f$. Thus to give an upper bound for $\|\partial_i(f - F)\|_{L^d}$, one term to be controlled is the $L^d$-norm of $\sum_j \mu_j \partial_i \omega_j$. But

(1.4) $$\left\|\sum_j |\partial_i \omega_j| |\mu_j|\right\|_{L^d} \leq C \left\|\sum_j 2^j \omega_j\right\|_{L^d};$$

it is therefore hopeless to conclude this way, since the right-hand side of (1.4) is even bigger than

$$\left\|\|2^j |\Delta_j f|\|_{\ell^1}\right\|_{L^d},$$

while one can only afford a bound by $\left\|\|2^j \Delta_j f\|_{\ell^2}\right\|_{L^d} \simeq \|\nabla f\|_{L^d}$. Bourgain and Brezis have a clever way out: if instead of $\left\|\sum_j 2^j \omega_j\right\|_{L^d}$ we only needed to estimate

$$\left\|\sum_j 2^j \omega_j \chi_{A_j}\right\|_{L^d}$$

where $A_j$ is the set defined by $A_j := \{x \in \mathbb{R}^d \colon \omega_j(x) > \sum_{t>0} 2^{-t} \omega_{j-t}(x)\}$ and $\chi_{A_j}$ is the characteristic function of the set $A_j$, then we would be in good shape because we have a pointwise bound

(1.5) $$\sum_j 2^j \omega_j \chi_{A_j} \leq 2 \sup_j 2^j \omega_j,$$

and the crucial estimate:

(1.6) $$\left\|\sup_{j \in \mathbb{Z}} (2^j \omega_j)\right\|_{L^p(\mathbb{R}^d)} \leq C 2^{\frac{\sigma(d-1)}{p}} \|\nabla f\|_{L^p(\mathbb{R}^d)}$$

for any $1 < p < \infty$. Thus they decompose

$$\Delta_j f(x) = \Delta_j f(x) \chi_{A_j^c}(x) + \Delta_j f(x) \chi_{A_j}(x) := g_j(x) + h_j(x)$$

so that $f = \sum_{j \in \mathbb{Z}} g_j + \sum_{j \in \mathbb{Z}} h_j$. They then proceed to approximate $g := \sum_{j \in \mathbb{Z}} g_j$ and $h := \sum_{j \in \mathbb{Z}} h_j$ by

$$\tilde{g} := \sum_{j \in \mathbb{Z}} g_j \prod_{j'>j} (1 - G_{j'}) \quad \text{and} \quad \tilde{h} := \sum_{j \in \mathbb{Z}} h_j \prod_{j'>j} (1 - U_{j'})$$

respectively, where $G_j$ and $U_j$ are some suitable controlling functions that satisfy pointwisely $g_j \leq G_j \leq 1$ and $h_j \leq U_j \leq 1$ (so that $\tilde{g}$ and $\tilde{h}$ are automatically bounded), whereas $G_j$ and $U_j$ are constructed from the $\omega_j$'s, so that the $L^d$-norms of $\partial_i(g - \tilde{g})$ and $\partial_i(h - \tilde{h})$ satisfy good estimates for $i = 1, \ldots, d-1$. Indeed, in [5], $h - \tilde{h}$ is written as a sum of products, which in turn allows a direct estimate of $\partial_i(h - \tilde{h})$ by the Leibniz rule; the heuristics centered around equations (1.5) and (1.6) suggest that $\|\partial_i(h - \tilde{h})\|_{L^d}$ may be small. On the other hand, $\partial_i(g - \tilde{g})$ is estimated using Littlewood-Paley inequalities, since it is a sum of pieces that are well-localized in frequency: indeed, note that

$$|\Delta_j f(x)|\chi_{A_j^c}(x) \leq \sum_{t>0} 2^{-t}\omega_{j-t}(x), \tag{1.7}$$

and while a derivative on the left hand side of (1.7) heuristically gains only $2^j$, a derivative on each term on the right hand side of (1.7) gains $2^{j-t}$, which is better when $t$ is large. It is this interplay that allows them to conclude with the estimate for $\partial_i(g - \tilde{g})$, and hence the proof of their theorem.

Now that we have recalled this basic strategy, we can address the difficulties we faced in extending the result of Bourgain and Brezis for $\dot{W}^{1,d}(\mathbb{R}^d)$, to the full Theorem 1.1 for $\dot{F}_q^{\alpha,p}(\mathbb{R}^d)$. The first difficulty arises when $\alpha > 1$: if we define the controlling functions $\omega_j$ as in [5] by using a sup-convolution, then the $\omega_j$ are at best Lipschitz, and in general may not be differentiated more than once. But an approximation theorem for $\dot{F}_q^{\alpha,p}(\mathbb{R}^d)$ naturally involves taking $\alpha$ derivatives, so a sup-convolution construction for the $\omega_j$'s cannot be expected to work when $\alpha > 1$. Following [33], where Bourgain and Brezis' result was extended to subelliptic settings, we overcome this by taking a discrete $\ell^p$ convolution instead; morally speaking, this means that we take

$$\omega_j(x) = \left[\sum_{r \in 2^{-j}\mathbb{Z}^d} \left(|\Delta_j f|(r) e^{-2^j|x''-r''|-2^{j-\sigma}|x'-r'|}\right)^p\right]^{1/p} \tag{1.8}$$

(here $r'$ and $r''$ are the first $\kappa$ and the last $d - \kappa$ variables of $r$ respectively, where $\kappa$ is defined as in Theorem 1.1; similarly for $x'$ and $x''$). For some technical reasons, this is not the precise definition of $\omega_j$ we will use; see (3.3) in Section 3 below for the precise construction of $\omega_j$. Once the correct definition of $\omega_j$ is in place, roughly speaking we would consider the sets

$$A_j := \{x \in \mathbb{R}^d \colon \omega_j(x) > \sum_{t>0} 2^{-\alpha t}\omega_{j-t}(x)\}$$

(note the dependence of this set on $\alpha$), and split

$$\Delta_j f(x) = \Delta_j f(x)\chi_{A_j^c}(x) + \Delta_j f(x)\chi_{A_j}(x) := g_j(x) + h_j(x)$$

as above (actually we would use a smooth version of $\chi_{A_j}$ instead of the sharp cut-off given by the characteristic function of $A_j$). We would then proceed as in [5] to approximate $\sum_{j \in \mathbb{Z}} h_j$ and $\sum_{j \in \mathbb{Z}} g_j$, except that several further difficulties must be overcome.

One of them is the proof of the analog of (1.6) in the case $q > p$. This arises, for instance, when we prove an approximation theorem for $\dot{W}^{k,d/k}(\mathbb{R}^d)$ with $d/2 < k < d$ (in which case $q = 2 > d/k = p$). In general, to prove Theorem 1.1 for $\dot{F}_q^{\alpha,p}(\mathbb{R}^d)$, we would like to prove an inequality of the form

$$\left\|\sup_{j \in \mathbb{Z}}(2^{\alpha j}\omega_j)\right\|_{L^p(\mathbb{R}^d)} \lesssim 2^{\frac{\sigma\kappa}{p}}\|f\|_{\dot{F}_q^{\alpha,p}(\mathbb{R}^d)}. \tag{1.9}$$

If $\omega_j$ was defined as in the putative definition (1.8), then morally speaking, the above inequality would admit an easy proof when $q \leq p$: indeed, heuristically we have

$$\omega_j(x) \simeq \left[\sum_{r \in \mathbb{Z}^d} \left(|\Delta_j f|(x - 2^{-j}r)e^{-|r''|-2^{-\sigma}|r'|}\right)^p\right]^{1/p},$$

so

$$\left\| \sup_{j \in \mathbb{Z}} (2^{\alpha j} \omega_j) \right\|_{L^p(\mathbb{R}^d)}^p = \int_{\mathbb{R}^d} \sup_{j \in \mathbb{Z}} (2^{\alpha j} \omega_j(x))^p dx$$

$$\leq \sum_{j \in \mathbb{Z}} \int_{\mathbb{R}^d} (2^{\alpha j} \omega_j(x))^p dx$$

$$= \sum_{j \in \mathbb{Z}} \sum_{r \in \mathbb{Z}^d} (e^{-|r''|-2^{-\sigma}|r'|})^p (2^{\alpha j})^p \int_{\mathbb{R}^d} |\Delta_j f(x - 2^{-j}r)|^p dx.$$

The last integral is equal to $\int_{\mathbb{R}^d} |\Delta_j f|(x)^p dx$, and $\sum_{r \in \mathbb{Z}^d} (e^{-|r''|-2^{-\sigma}|r'|})^p \lesssim 2^{\sigma \kappa}$. Thus

$$\left\| \sup_{j \in \mathbb{Z}} (2^{\alpha j} \omega_j) \right\|_{L^p(\mathbb{R}^d)}^p \lesssim 2^{\sigma \kappa} \int_{\mathbb{R}^d} \sum_{j \in \mathbb{Z}} (2^{\alpha j} |\Delta_j f(x)|)^p dx$$

$$\lesssim 2^{\sigma \kappa} \int_{\mathbb{R}^d} \left( \sum_{j \in \mathbb{Z}} (2^{\alpha j} |\Delta_j f(x)|)^q \right)^{p/q} dx$$

where in the last line we have used the embedding $\ell^q \hookrightarrow \ell^p$ if $q \leq p$. This would prove (1.9) when $q \leq p$, under the putative definition (1.8) of $\omega_j$. Unfortunately this simple argument is insufficient in handling the case when $q > p$. We found a way out using a logarithmic bound for some vector-valued 'shifted' maximal functions (see Corollary 10.3), which we prove using an old argument going back to Zó ([35]). We then get a slightly weaker bound than (1.9), one that is off by a logarithmic factor (see Proposition 4.7), but that is still sufficient for our purpose.

A second difficulty arises when $\alpha$ is not an integer or when $q \neq 2$. Recall that in one step, Bourgain and Brezis estimated $\partial_i(h - \tilde{h})$ in $L^d(\mathbb{R}^d)$ by writing it as a sum of products, and then using the ordinary Leibniz rule. In our case, we need to estimate $\partial_i(h - \tilde{h})$ in $\dot{F}_q^{\alpha-1,p}(\mathbb{R}^d)$, which is defined only via Littlewood-Paley projections when $\alpha$ is not an integer, or when $q \neq 2$. Thus we must know how to estimate the derivative of a sum of products within the realm of Littlewood-Paley theory. If it were not for the sum involved, we could just apply the fractional Leibniz rule for the space $\dot{F}_q^{\alpha-1,p}(\mathbb{R}^d)$. But since the sum is present, we found it easier to proceed directly, without resorting to the fractional Leibniz rule. It may also be worth noting here that we run into an additional difficulty, in the case $0 < \alpha < 1$: we find it necessary then to exploit some additional cancellations offered by the Littlewood-Paley projections $\Delta_j$'s, when we deal with certain high frequency components of $h - \tilde{h}$ (see the introduction of the parameter $T$ in Section 5 when $0 < \alpha < 1$).

A final difficulty arises when $\alpha \in (0, 1/2]$. In this case, $\alpha$ is rather small, so the set $A_j^c$, given by $A_j^c = \{x \in \mathbb{R}^d : \omega_j(x) \leq \sum_{t>0} 2^{-\alpha t} \omega_{j-t}(x)\}$, is relatively large. As a result, $g_j := \Delta_j f \cdot \chi_{A_j^c}$ is relatively large, and one expects it to be relatively harder to estimate $\partial_i(g - \tilde{g})$ in $\dot{F}_q^{\alpha-1,p}(\mathbb{R}^d)$. This is manifested in our need to introduce a parameter $a_\alpha$ in Section 6 (see Proposition 6.1), which is smaller than 1 when $\alpha \in (0, 1/2]$.

Let us end up this introduction with three open problems:

**Open problem 1.3.** *The condition $\kappa < \min(p, d)$ in the statements of Theorems* 1.1 *and* 1.2 *may not be necessary in general. In* [17], *Maz'ya proves that for every vector function $X \in$* $\dot{F}_2^{d/2,2}(\mathbb{R}^d; \mathbb{R}^d)$, *there exists* $Y \in (\dot{F}_2^{d/2,2} \cap L^\infty)(\mathbb{R}^d; \mathbb{R}^d)$ *and a scalar function* $u \in \dot{F}_2^{1+d/2,2}(\mathbb{R}^d)$ *such that*

$$X = Y + \nabla u.$$

*This coincides with the statement of Theorem 1.2 when $p = q = 2$, $\alpha = d/2$ and $l = 1$, except that this set of parameters is not covered by our assumptions when $d \geq 3$. Indeed, the condition $l \in [\![d-\kappa, d-1]\!]$ cannot be satisfied in that case: this would require $\kappa = d-1$, which is impossible in view of the conditions $\kappa < p = 2$ and $d \geq 3$. Is it true that Theorems 1.1 and 1.2 remain true when the condition $\kappa < \min\{p, d\}$ is replaced by $\kappa < d$?*

**Open problem 1.4.** *In [33], the conclusion of [5] is extended to a subelliptic context, namely the case of the Heisenberg groups endowed with a subelliptic Laplacian. The extension of Theorems 1.1 and 1.2 to the case of the Heisenberg group is an open problem.*

**Open problem 1.5.** *It is likely that Theorem 1.2 can be extended to the case of smooth bounded domains in $\mathbb{R}^d$, in the spirit of [6].*

The paper is organized as follows. After gathering instrumental facts about Triebel-Lizorkin spaces and maximal functions in Section 2, we describe the approximating function $F$ in Theorem 1.1 in Section 3. Section 4 is devoted to proving key estimates for the $\omega_j$'s, which are then used to derive bounds for $h - \tilde{h}$ (resp. $g - \tilde{g}$) in Section 5 (resp. Section 6). The proof of Theorem 1.1 is completed in Section 7, while Theorem 1.2 is established in Section 8.

Throughout the paper, if two quantities $A(f)$ and $B(f)$ depend on a function $f$ ranging over some space $L$, the notation $A(f) \lesssim B(f)$ means that there exists $C > 0$ such that $A(f) \leq CB(f)$ for all $f \in L$, while $A(f) \simeq B(f)$ means that $A(f) \lesssim B(f) \lesssim A(f)$. The Euclidean ball centered at $0$ with radius $r$ will be denoted $B_r$.

**Acknowledgment.** Wang was partially supported by NSF Grant No. DMS-1612015. Yung was partially supported by the Early Career Grant CUHK24300915 from the Hong Kong Research Grant Council.

## 2. Preliminaries

For a brief overview on homogeneous Triebel-Lizorkin spaces, we refer to [25, Chapter 5], [3] and also [21, Chapter 2].

2.1. **The Triebel-Lizorkin spaces.** We fix a function $\Delta \in \mathcal{S}(\mathbb{R}^d)$ such that[1]

(2.1) $$\hat{\Delta} \in C_c^\infty(B_2 \setminus B_{\frac{1}{2}}),$$

and

(2.2) $$\sum_{j \in \mathbb{Z}} \hat{\Delta}(2^j x) = 1, \qquad \forall x \in \mathbb{R}^d \setminus \{0\}.$$

Notice that assumption (2.1) yields that, for every polynomial $P$,

(2.3) $$\int_{\mathbb{R}^d} P(x) \Delta(x)\, dx = 0.$$

For all $j \in \mathbb{Z}$ and all $x \in \mathbb{R}^d$, define $\Delta_j(x) := 2^{jd} \Delta(2^j x)$.

---

[1]The function $\Delta$ can be obtained as follows. Let $\rho \in C_c^\infty(B_2 \setminus B_{\frac{1}{2}})$ such that $\rho \equiv 1$ on $C_c^\infty(B_{\frac{3}{2}} \setminus B_{\frac{3}{4}})$ and $0 \leq \rho \leq 1$ on $\mathbb{R}^d$. Then $1 \leq \sum_{j \in \mathbb{Z}} \rho(2^j x) \leq 2$ on $\mathbb{R}^d \setminus \{0\}$. We then define $\hat{\Delta}(x) = \frac{\rho(x)}{\sum_{j \in \mathbb{Z}} \rho(2^j x)}$ on $\mathbb{R}^d \setminus \{0\}$ and $\hat{\Delta}(0) = 0$.

Let
$$\mathcal{Z}(\mathbb{R}^d) = \{f \in \mathcal{S}(\mathbb{R}^d) : \partial^\gamma \hat{f}(0) = 0, \forall \gamma \in \mathbb{N}^d\}.$$
The dual space $\mathcal{Z}'(\mathbb{R}^d)$ of this closed subspace of $\mathcal{S}(\mathbb{R}^d)$ can be identified as the set $\{f|_\mathcal{Z}, f \in \mathcal{S}'(\mathbb{R}^d)\}$, or equivalently as the factor space $\mathcal{S}'(\mathbb{R}^d)/\mathcal{P}(\mathbb{R}^d)$, where $\mathcal{P}(\mathbb{R}^d)$ is the collection of all polynomials on $\mathbb{R}^d$.

For all $f \in \mathcal{Z}'(\mathbb{R}^d)$, let $\Delta_j f := f * \Delta_j$; this is well-defined for $f \in \mathcal{Z}'(\mathbb{R}^d) = \mathcal{S}'(\mathbb{R}^d)/\mathcal{P}(\mathbb{R}^d)$ since the Fourier transform of a polynomial is supported in $\{0\}$. Moreover, it is a straightforward consequence of the Paley-Wiener theorem that $\Delta_j f$ belongs to $L^p(\mathbb{R}^d)$ for all $p \in [1, \infty]$.

**Definition 2.1.** Let $\alpha \in \mathbb{R}$ and $p, q \in (1, \infty)$. Let $f \in \mathcal{Z}'(\mathbb{R}^d)$. Say that $f \in \dot{F}_q^{\alpha,p}(\mathbb{R}^d)$ (or $\dot{F}_q^{\alpha,p}$) if and only if
$$\|f\|_{\dot{F}_q^{\alpha,p}(\mathbb{R}^d)} := \|\|2^{\alpha j}\Delta_j f\|_{\ell^q(\mathbb{Z})}\|_{L^p(\mathbb{R}^d)} < \infty.$$

In view of (2.1) and (2.2), we have

(2.4) $$\forall k \in \mathbb{Z}, \forall f \in \mathcal{Z}', \quad \Delta_k f = \sum_{j \in \mathbb{Z}} \Delta_k \Delta_j f = \sum_{|j-k| \le 1} \Delta_k \Delta_j f.$$

This implies

**Proposition 2.2.** Let $\alpha > 0$ and $1 < p, q < \infty$. Then for all $f \in \dot{F}_q^{\alpha,p}(\mathbb{R}^d)$,

(2.5) $$f = \sum_{j=-\infty}^{\infty} \Delta_j f$$

where the series converges in $\dot{F}_q^{\alpha,p}$.

Another useful property is given by the following proposition:

**Proposition 2.3.** [3, Proposition 5] For every $\alpha, p, q \in \mathbb{R}$, for every $f \in \mathcal{Z}'(\mathbb{R}^d)$,
$$\|f\|_{\dot{F}_q^{\alpha,p}} \simeq \sum_{i=1}^{d} \|\partial_i f\|_{\dot{F}_q^{\alpha-1,p}}.$$

**Remark 2.4.** We only need to prove Theorem 1.1 under the additional assumption that $f \in \dot{F}_q^{\alpha,p}$ has only finitely many $\Delta_j f$ different from 0.

Indeed, assume that the theorem is true for such distributions $f$. Then for an arbitrary $f \in \dot{F}_q^{\alpha,p}$, we consider for every $J \in \mathbb{N}$ the distribution
$$f_J := \sum_{j=-J}^{J} \Delta_j f.$$

For every $\delta > 0$, we thus get a function $F_J \in \dot{F}_q^{\alpha,p} \cap L^\infty$ which satisfies

(2.6) $$\sum_{i=1}^{\kappa} \|\partial_i(f_J - F_J)\|_{\dot{F}_q^{\alpha-1,p}} \le \delta \|f_J\|_{\dot{F}_q^{\alpha,p}},$$

and

(2.7) $$\|F_J\|_{L^\infty} + \|F_J\|_{\dot{F}_q^{\alpha,p}} \le C_\delta \|f_J\|_{\dot{F}_q^{\alpha,p}}.$$

Proposition 2.2 implies that the sequence $(f_J)_{J \in \mathbb{N}}$ strongly converges to $f$ in $\dot{F}_q^{\alpha,p}$. Hence, the sequence $(F_J)_{J \in \mathbb{N}}$ is bounded in $\dot{F}_q^{\alpha,p} \cap L^\infty$. We can extract a subsequence (still denoted by $F_j$) which converges to some $F$ weakly* in $L^\infty$, and thus also in $\mathcal{Z}'$. By the Fatou property [3, Proposition 7], $F \in L^\infty \cap \dot{F}_q^{\alpha,p}$ and (2.7) remains true with $F$ and $f$ instead of $F_J$ and $f_J$. Since $\partial_i(f_J - F_J)$ also converges weakly* in $\mathcal{Z}'$, the Fatou property again implies that (2.6) remains true for $f$ and $F$.

We assume henceforth in all the sequel of the paper that $f$ is such that only finitely many $\Delta_j f$ are different from 0.

2.2. **Inequalities involving the Hardy-Littlewood maximal function.** For all functions $g \in L_{loc}^1(\mathbb{R}^d)$ and all $x \in \mathbb{R}^d$, define the Hardy-Littlewood functional by

$$\mathcal{M}g(x) := \sup_{B \ni x} \frac{1}{|B|} \int_B |g(y)|\, dy,$$

where the supremum is taken over all Euclidean balls of $\mathbb{R}^d$ containing $x$. Let us summarize the properties of $\mathcal{M}$ which will be used in the sequel:

**Proposition 2.5.** *The Hardy-Littlewood functional satisfies the following properties:*

(1) $\mathcal{M}$ *is of weak type* $(1,1)$ *and* $L^p$-*bounded for all* $p \in (1, \infty]$, [2]
(2) *one also has the vector-valued version of the previous assertion: for all* $p, q \in (1, \infty)$, *for all* $(g_j)_{j \in \mathbb{Z}} \in L^p(\mathbb{R}^d; \ell^q(\mathbb{Z}))$,

(2.8) $$\|\|\mathcal{M}g_j\|_{\ell^q(j)}\|_{L^p} \lesssim \|\|g_j\|_{\ell^q(j)}\|_{L^p},$$

(3) *for all* $p \in [1, \infty]$, *all functions* $g \in L^p(\mathbb{R}^d)$, *all decreasing functions* $\varphi : [0, \infty) \to [0, \infty)$ *such that* $A := \int_{\mathbb{R}^d} \varphi(\|y\|) dy < \infty$,[3] *and all measurable functions* $\phi$ *such that* $|\phi(y)| \leq \varphi(\|y\|)$ *for all* $y \in \mathbb{R}^d$, *the convolution* $g * \phi$ *is defined almost everywhere and one has*

(2.9) $$|g * \phi(x)| \lesssim A\mathcal{M}g(x).$$

*Proof.* See ([23, Chapter 1, Theorem 1], [23, Chapter 2, Theorem 1]) and [22, Chapter 3, Theorem 2(a)]) respectively. □

**Remark 2.6.** Note that (2.9) applies in particular when $\phi \in \mathcal{S}(\mathbb{R}^d)$, since every Schwartz function on $\mathbb{R}^d$ can be dominated by a radially decreasing integrable function.

**Proposition 2.7.** *For all* $\gamma \in \mathbb{N}^d$, $j \in \mathbb{Z}$ *and* $x \in \mathbb{R}^d$,

(2.10) $$|\partial^\gamma \Delta_j f(x)| \lesssim 2^{|\gamma|j} \mathcal{M}\Delta_j f(x).$$

*Moreover, for all* $\alpha > 0$ *and* $1 < p, q < \infty$ *such that* $\alpha p = d$,

(2.11) $$|\partial^\gamma \Delta_j f(x)| \lesssim 2^{|\gamma|j} \|f\|_{\dot{F}_q^{\alpha,p}}.$$

*The implicit constants in both inequalities do not depend on* $x$, $\gamma$, $j$ *nor on* $f$.

---

[2] Note that $\mathcal{M}$ is a sublinear operator. That $\mathcal{M}$ is $L^p$-bounded (resp. is of weak type $(1,1)$) means that $\|\mathcal{M}g\|_p \lesssim \|g\|_p$ (resp. that, for all $\lambda > 0$, $|\{\mathcal{M}g > \lambda\}| \lesssim \frac{1}{\lambda}\|g\|_1$).

[3] Here and after, $\|\cdot\|$ stands for the Euclidean norm.

As a particular case of (2.11) where we take $\gamma = 0$, we obtain the following Bernstein inequality when $\alpha p = d$:

$$\|\Delta_j f\|_{L^\infty} \lesssim \|f\|_{\dot{F}_q^{\alpha,p}}. \tag{2.12}$$

The implicit constant in (2.12) only depends on $\alpha$ and $p$ (but neither on $j$ nor on $f$).

*Proof of Proposition* 2.7. In view of (2.4), we have

$$|\partial^\gamma \Delta_j f(x)| \lesssim \sum_{|j-k|\leq 1} 2^{|\gamma|k} |(\partial^\gamma \Delta)_k \Delta_j f| \lesssim 2^{|\gamma|j} \sum_{|j-k|\leq 1} |(\partial^\gamma \Delta)_k \Delta_j f|,$$

where $(\partial^\gamma \Delta)_k(x) = 2^{kd}(\partial^\gamma \Delta)(2^k x)$. Taking (2.9) into account and applying Remark 2.6 to $\phi = (\partial^\gamma \Delta)_k$, this yields

$$|\partial^\gamma \Delta_j f(x)| \lesssim 2^{|\gamma|j} \mathcal{M} \Delta_j f(x).$$

This proves the first assertion. It follows therefrom that

$$|\partial^\gamma \Delta_j f(x)| \lesssim 2^{|\gamma|j} \|\Delta_j f\|_{L^\infty}. \tag{2.13}$$

Using (2.4) again, we have

$$\|\Delta_j f\|_{L^\infty} \leq \sum_{|j-k|\leq 1} \|\Delta_k \Delta_j f\|_{L^\infty}.$$

Hölder's inequality then implies

$$\|\Delta_j f\|_{L^\infty} \leq \sum_{|k-j|\leq 1} \|\Delta_k\|_{L^{p'}} \|\Delta_j f\|_{L^p}.$$

Using a change of variables and the expression of $\Delta_k$, we thus get

$$\|\Delta_j f\|_{L^\infty} \lesssim 2^{\frac{jd}{p}} \|\Delta_j f\|_{L^p} = 2^{\alpha j} \|\Delta_j f\|_{L^p} \leq \|f\|_{\dot{F}_q^{\alpha,p}}, \tag{2.14}$$

where we have used that $\alpha p = d$ and the definition of $\dot{F}_q^{\alpha,p}$. Inequality (2.11) is now a consequence of (2.13) and (2.14).

$\square$

## 3. The approximations of $f$

The present section is devoted to the definition of the function $F$ in Theorem 1.1. Let $\alpha > 0$ and $p > 1$ such that $\alpha p = d$. Let $f \in \dot{F}_q^{\alpha,p}(\mathbb{R}^d)$, $\delta > 0$ and $\sigma$ be a large positive integer to be chosen (only depending on $\delta$). For $x = (x_1, \ldots, x_d)$, we define $x_\sigma := (2^{-\sigma} x_1, \ldots, 2^{-\sigma} x_\kappa, x_{\kappa+1}, \ldots, x_d)$. The parameter $\sigma$ discriminates the *good* directions $x_1, \ldots, x_\kappa$ from the other ones. In particular, when one differentiates a function of the form $x \mapsto u(x_\sigma)$ along a good direction, an additional factor $2^{-\sigma}$ arises.

Let $E$ be the Schwartz function defined by

$$E(x) := e^{-(1+\|x_\sigma\|^2)^{\frac{1}{2}}}. \tag{3.1}$$

For all $j \in \mathbb{Z}$ and all $x \in \mathbb{R}^d$, let $E_j(x) := 2^{jd} E(2^j x)$.

Define also

$$T(x) := \min\left(1, \|x\|^{-(d+1)}\right) \tag{3.2}$$

and

$$T_j(x) := 2^{jd} T(2^j x)$$

for all $j \in \mathbb{Z}$ and all $x \in \mathbb{R}^d$.

We introduce an auxiliary function which can be seen as a substitute of $|\Delta_j f|$, $j \in \mathbb{Z}$:

$$\omega_j(x) := \left( \sum_{r \in \mathbb{Z}^d} [T_j |\Delta_j f|(2^{-j}r) E(2^j x - r)]^p \right)^{\frac{1}{p}}, \quad x \in \mathbb{R}^d. \tag{3.3}$$

Here and in the sequel, we use the notation $T_j |\Delta_j f|$ for the convolution $T_j * |\Delta_j f|$. We will prove that $\omega_j$ inherits the $L^\infty$ bounds of $|\Delta_j f|$. More precisely,

$$\|\Delta_j f\|_{L^\infty} \lesssim \|\omega_j\|_{L^\infty} \lesssim 2^{\kappa\sigma} \|\Delta_j f\|_{L^\infty}.$$

In contrast to $|\Delta_j f|$, $\omega_j$ is smooth, as a discrete $\ell^p$ convolution. Moreover, it behaves differently with respect to good and bad coordinates. This allows to obtain improved estimates on its derivatives along good directions.

**Remark 3.1.** Notice that, if, for some $x \in \mathbb{R}^d$ and some $j \in \mathbb{Z}$, $\omega_j(x) = 0$, then the definition of $\omega_j$ yields that $T_j |\Delta_j f|(2^{-j}r) = 0$ for all $r \in \mathbb{Z}$. Since $T_j$ is positive everywhere, it follows that $\Delta_j f$ has to vanish on all $\mathbb{R}^d$, which entails that $\omega_j(y) = 0$ for all $y \in \mathbb{R}^d$.

Let $R >> \sigma$ be another positive integer to be chosen. Let us consider a smooth function $\zeta_j$ approximating the characteristic function of the set

$$\left\{ x \in \mathbb{R}^d;\ 2^{\alpha j} \omega_j(x) \leq \frac{1}{2} \sum_{k<j,\, k\equiv j \,(\text{mod } R)} 2^{\alpha k} \omega_k(x) \right\}.$$

More specifically, notice first that, if the function $\sum_{k<j, k\equiv j (\text{mod } R)} 2^{\alpha k} \omega_k$ vanishes at some point $x \in \mathbb{R}^d$, then it identically vanishes (see Remark 3.1). In order to define $\zeta_j$, we thus fix a smooth function $\zeta : [0,\infty) \to [0,1]$ such that $\zeta \equiv 1$ on $[0, \frac{1}{2}]$ and $\zeta \equiv 0$ on $[1,\infty)$ and we define, for all $j \in \mathbb{Z}$,

$$\zeta_j := \begin{cases} \zeta \left( \dfrac{2^{\alpha j} \omega_j}{\sum_{k<j, k\equiv j (\text{mod } R)} 2^{\alpha k} \omega_k} \right) & \text{if } \sum_{k<j, k\equiv j (\text{mod } R)} 2^{\alpha k} \omega_k \not\equiv 0, \\ 0 & \text{otherwise}. \end{cases}$$

We split $f = \sum_j \Delta_j f$ into the sum of two functions: $f = g + h$, given by

$$h := \sum_{j=-\infty}^{\infty} h_j, \quad g := \sum_{j=-\infty}^{\infty} g_j, \tag{3.4}$$

with

$$h_j(x) := (1 - \zeta_j(x)) \Delta_j f(x), \quad g_j(x) := \zeta_j(x) \Delta_j f(x).$$

The approximating function $F$ in Theorem 1.1 is also defined as the sum of two functions: $F = \tilde{g} + \tilde{h}$, where

$$\tilde{h} := \sum_{j=-\infty}^{\infty} h_j \prod_{j'>j} (1 - U_{j'}) \quad \text{with} \quad U_j = (1 - \zeta_j) \omega_j \tag{3.5}$$

and

$$\tilde{g} = \sum_{c=0}^{R-1} \sum_{j \equiv c (\mathrm{mod}\, R)} g_j \prod_{\substack{j' > j \\ j' \equiv c (\mathrm{mod}\, R)}} (1 - G_{j'}) \quad \text{with} \quad G_j := \sum_{\substack{t > 0 \\ t \equiv 0 (\mathrm{mod}\, R)}} 2^{-\alpha t} \omega_{j-t}. \tag{3.6}$$

The definition of $\tilde{g}$ involves some infinite products, the convergence of which will be discussed at the end of Section 4, while, as we shall see, the products involved in the definition of $\tilde{h}$ are actually finite. We will show that $F$ satisfies all the conclusions of Theorem 1.1, provided that $\|f\|_{\dot{F}_q^{\alpha,p}}$ is sufficiently small. The latter can be assumed without loss of generality, as explained at the end of this section, see Proposition 3.3 below.

The definitions of $\tilde{h}$ and $\tilde{g}$ are inspired by [5] and [33]. They are motivated by two crucial facts: the Bernstein inequality (2.12) and the following algebraic identity, see e.g. [33, Proposition 6.1]:

**Lemma 3.2.** *Let $(a_k)_{k \geq \mathbb{Z}}$ be a sequence of complex numbers. Assume that, for some integer $k_0 \in \mathbb{Z}$, $a_k = 0$ whenever $|k| > k_0$. Then, for all $j \in \mathbb{Z}$,*

$$1 = \sum_{j' > j} a_{j'} \prod_{j < j'' < j'} (1 - a_{j''}) + \prod_{j' > j} (1 - a_{j'}).$$

In particular, the above identity implies that if $0 \leq a_j \leq 1$ for every $j \in \mathbb{Z}$, then

$$\left| \sum_{j' > j} a_{j'} \prod_{j < j'' < j'} (1 - a_{j''}) \right| \leq 1. \tag{3.7}$$

The functions $U_j$ in the definition of $\tilde{h}$ are constructed in order to satisfy

$$\|h_j\|_{L^\infty} \lesssim \|U_j\|_{L^\infty} \lesssim \|\omega_j\|_{L^\infty} \lesssim 2^{\kappa \sigma} \|\Delta_j f\|_{L^\infty}.$$

Taking $a_j = U_j$ in (3.7) and using Bernstein inequality, one can see that when $\|f\|_{\dot{F}_q^{\alpha,p}}$ is sufficiently small, $0 \leq U_j \leq 1$ and thus $\|\tilde{h}\|_{L^\infty} \lesssim 1$. A similar computation can be made with $\tilde{g}$. This will imply the desired $L^\infty$ estimate on $\tilde{F}$.

Regarding the $\dot{F}_q^{\alpha,p}$ estimates, the strategies for $\tilde{h}$ and $\tilde{g}$ follow two different paths. For every $x \in \mathbb{R}^d$, $h(x)$ is the sum of the largest Littlewood-Paley projections $|\Delta_j f(x)|$. Roughly speaking, this is exploited to reduce the sum of these projections to only one term. More specifically, we will use the following fact:

$$\sum_{m \in \mathbb{Z}} 2^{\alpha m} \omega_m \chi_{2^{\alpha m} \omega_m > \frac{1}{2} \sum_{k < m, k \equiv m (\mathrm{mod}\, R)} 2^{\alpha k} \omega_k} \leq 3R \sup_{m \in \mathbb{Z}} 2^{\alpha m} \omega_m. \tag{3.8}$$

Indeed, writing

$$\sum_{m \in \mathbb{Z}} 2^{\alpha m} \omega_m \chi_{2^{\alpha m} \omega_m > \frac{1}{2} \sum_{k < m, k \equiv m (\mathrm{mod}\, R)} 2^{\alpha k} \omega_k} = \sum_{j=0}^{R-1} \sum_{m \equiv j (\mathrm{mod}\, R)} 2^{\alpha m} \omega_m \chi_{2^{\alpha m} \omega_m > \frac{1}{2} \sum_{k < m, k \equiv j (\mathrm{mod}\, R)} 2^{\alpha k} \omega_k},$$

we consider for every $j = 0, \ldots, R-1$ the largest index $m_j$ in the sum $\sum_{m \equiv j (\mathrm{mod}\, R)} \cdots$ above such that the corresponding term $2^{\alpha m} \omega_m \chi_{2^{\alpha m} \omega_m > \frac{1}{2} \sum_{k < m, k \equiv j (\mathrm{mod}\, R)} 2^{\alpha k} \omega_k}$ is $> 0$ (such an index $m_j$ exists since we have assumed that only finitely many $\Delta_k f$ are different from 0). Then

$$\sum_{m \equiv j (\mathrm{mod}\ R)} 2^{\alpha m} \omega_m \chi_{2^{\alpha m} \omega_m > \frac{1}{2} \sum_{k<m, k \equiv m (\mathrm{mod}\ R)} 2^{\alpha k} \omega_k} \leq 2^{\alpha m_j} \omega_{m_j} + \sum_{k < m_j, k \equiv j (\mathrm{mod}\ R)} 2^{\alpha k} \omega_k$$

$$\leq 3 \cdot 2^{\alpha m_j} \omega_{m_j} \leq 3 \sup_{m \in \mathbb{Z}} 2^{\alpha m} \omega_m$$

from which (3.8) follows. The estimate of the $\dot{F}_q^{\alpha,p}$ norm of the right hand side of (3.8) is the most delicate part in the $\dot{F}_q^{\alpha,p}$ approximation of $h$ by $\tilde{h}$. This is the object of Proposition 4.7 below. Let us also mention that the good derivatives play a central role in this first part of the approximation.

The $\dot{F}_q^{\alpha,p}$ estimate of $g - \tilde{g}$ is less elaborate. As explained in the introduction, it is obtained using Littlewood-Paley inequalities. Here, the role of $R$ becomes crucial.

In order to carry out the above arguments rigorously, we need to assume that $\|f\|_{\dot{F}_q^{\alpha,p}}$ is sufficiently small. This is not a restriction since Theorem 1.1 is a consequence of the following (apparently weaker) statement:

**Proposition 3.3.** *Let $\alpha > 0$ and $p, q \in (1, \infty)$ such that $\alpha p = d$. Let $\kappa$ be the largest positive integer that satisfies $\kappa < \min\{p, d\}$. Then for every $\delta > 0$, there exists $\eta_\delta > 0$ such that for every $f \in \dot{F}_q^{\alpha,p}(\mathbb{R}^d)$ with $\|f\|_{\dot{F}_q^{\alpha,p}} \leq \eta_\delta$, there exist $F \in \dot{F}_q^{\alpha,p} \cap L^\infty(\mathbb{R}^d)$ and a constant $D_\delta > 0$ with $D_\delta$ independent of $f$, such that*

$$(3.9) \qquad \sum_{i=1}^{\kappa} \|\partial_i (f - F)\|_{\dot{F}_q^{\alpha-1,p}} \leq \delta \|f\|_{\dot{F}_q^{\alpha,p}} + D_\delta \|f\|_{\dot{F}_q^{\alpha,p}}^2,$$

*and*

$$(3.10) \qquad \|F\|_{L^\infty} + \|F\|_{\dot{F}_q^{\alpha,p}} \leq D_\delta.$$

We proceed to explain how Proposition 3.3 implies Theorem 1.1. Let $\delta > 0$. By Proposition 3.3, there exist $\eta_\delta > 0$ and $D_\delta > 0$ satisfying the above properties.

Let $f \in \dot{F}_q^{\alpha,p}$, $f \not\equiv 0$. We then apply Proposition 3.3 to the function

$$\tilde{f} := \frac{\min\left(\eta_\delta, \frac{\delta}{D_\delta}\right)}{\|f\|_{\dot{F}_q^{\alpha,p}}} f.$$

We thus obtain a function $\tilde{F}$ which satisfies (3.9) and (3.10), with $\tilde{f}$ instead of $f$. Finally, we set

$$F := \frac{\|f\|_{\dot{F}_q^{\alpha,p}}}{\min\left(\eta_\delta, \frac{\delta}{D_\delta}\right)} \tilde{F}.$$

Then multiplying the estimates by $\|f\|_{\dot{F}_q^{\alpha,p}} / \min\left(\eta_\delta, \frac{\delta}{D_\delta}\right)$ and using that $\left\|\tilde{f}\right\|_{\dot{F}_q^{\alpha,p}} = \min\left(\eta_\delta, \frac{\delta}{D_\delta}\right)$ yields

$$\sum_{i=1}^{\kappa} \|\partial_i (f - F)\|_{\dot{F}_q^{\alpha-1,p}} \leq \delta \|f\|_{\dot{F}_q^{\alpha,p}} + D_\delta \min\left(\eta_\delta, \frac{\delta}{D_\delta}\right) \|f\|_{\dot{F}_q^{\alpha,p}} \leq 2\delta \|f\|_{\dot{F}_q^{\alpha,p}},$$

and

$$\|F\|_{L^\infty} + \|F\|_{\dot{F}_q^{\alpha,p}} \leq \frac{D_\delta}{\min\left(\eta_\delta, \frac{\delta}{D_\delta}\right)} \|f\|_{\dot{F}_q^{\alpha,p}}.$$

This proves Theorem 1.1 with $C_\delta = D_{\delta/2}/\min\left(\eta_{\delta/2}, \frac{\delta}{2D_{\delta/2}}\right)$.

A word about notations is in order. In the above, we have defined the functions $E_j$, $T_j$, $\Delta_j$, $\omega_j$, $\zeta_j$, $h_j$, $g_j$, $U_j$ and $G_j$. Morally speaking, all these are localized in frequency to $|\xi| \lesssim 2^j$. Some like $E_j, T_j$ and $\Delta_j$ are $L^1$-renormalized dilations of a fixed function (in particular, we note in passing that they satisfy

$$\|E_j\|_{L^p} = 2^{\frac{jd}{p'}} \|E\|_{L^p}$$

for all $p$; similarly for $T_j$ and $\Delta_j$). The others are not dilations of a fixed function, but if we take $k$ derivatives of $\omega_j$, $\zeta_j$, $h_j$, $g_j$, $U_j$ or $G_j$, we can obtain an upper bound that involves a factor $2^{jk}$. This will be made explicit in the next three sections, which are devoted to the proof of Proposition 3.3.

## 4. Properties of $\omega_j$

In this section, we collect all the estimates of $\omega_j$ needed in the sequel.

### 4.1. Pointwise estimates.
First we have (see [33, Section 9]):

**Proposition 4.1.** *For all $j \in \mathbb{Z}$ and all $x \in \mathbb{R}^d$,*

$$(4.1) \qquad \omega_j(x) \simeq \left(\sum_{r \in \mathbb{Z}^d} \left(T_j|\Delta_j f|(x + 2^{-j}r)E(r)\right)^p\right)^{1/p}$$

$$(4.2) \qquad |\Delta_j f|(x) \lesssim \omega_j(x).$$

The proof of Proposition 4.1 relies on the following estimate for the function $T$ defined in (3.2):

**Lemma 4.2.** *For every $x, y \in \mathbb{R}^d$, if $\|y\| \leq \sqrt{d}$, then $T(x + y) \lesssim T(x)$.*

*Proof of Lemma 4.2:* Note that $T(x) \geq (2\sqrt{d})^{-(d+1)}$ for all $x \in \mathbb{R}^d$ with $\|x\| \leq 2\sqrt{d}$. Thus

$$T(x+y) \leq 1 \leq (2\sqrt{d})^{d+1} T(x)$$

whenever $\|x\| \leq 2\sqrt{d}$. On the other hand, if $\|x\| > 2\sqrt{d}$ and $\|y\| \leq \sqrt{d}$, then $\|x+y\| \geq \|x\|/2$, so

$$T(x+y) = \|x+y\|^{-(d+1)} \leq 2^{d+1} \|x\|^{-(d+1)} = 2^{d+1} T(x).$$

This shows $T(x+y) \leq C T(x)$ with $C = (2\sqrt{d})^{d+1}$. □

*Proof of Proposition 4.1:* The proof of (4.1) is analogous to the one of [33, Proposition 9.2]. The only difference is that the function $S_{j+N}$ introduced in [33] is now replaced by the function $T_j$ which satisfies the same (crucial) property as $S_{j+N}$, namely: for every $x, y \in \mathbb{R}^d$, if $\|y\| \leq 2^{-j}\sqrt{d}$,

$$(4.3) \qquad T_j(x+y) \simeq T_j(x).$$

In turn, this follows from the definition of $T_j$ and Lemma 4.2. From (4.3) one deduces that

$$(4.4) \qquad T_j|\Delta_j f|(x - 2^{-j}(y+r)) \simeq T_j|\Delta_j f|(x - 2^{-j}r)$$

whenever $x \in \mathbb{R}^d$, $r \in \mathbb{Z}^d$ and $\|y\| \leq \sqrt{d}$. Arguing as in the proof of [33, Proposition 9.2], we rewrite $\omega_j$ as

$$\omega_j(x) = \left(\sum_{r \in \mathbb{Z}^d} (T_j|\Delta_j f|(x - 2^{-j}(y+r))E(y+r))^p\right)^{1/p}$$

where $y \in [0,1)^d$ is the 'fractional part' of $2^j x$. In particular, $\|y\| \leq \sqrt{d}$. The estimate (4.1) then follows from (4.4) and the fact that whenever $r \in \mathbb{Z}^d$,

$$E(y+r) \simeq E(r). \tag{4.5}$$

Let now $K$ be a Schwartz function on $\mathbb{R}^d$, whose Fourier transform $\hat{K}$ is identically 1 on $B_2$, and vanishes outside $B_3$. Then by (2.1), for every $\xi \in \mathbb{R}^d$,

$$\hat{\Delta}(\xi)\hat{K}(\xi) = \hat{\Delta}(\xi).$$

Hence, $\Delta_j f = \Delta_j f * K_j$ where $K_j(x) = 2^{jd} K(2^j x)$ for all $x \in \mathbb{R}^d$. Moreover, since $K \in \mathcal{S}$, there exists $C > 0$ such that, for all $x \in \mathbb{R}^d$, $|K(x)| \leq CT(x)$. We deduce therefrom $|K_j| \leq CT_j$ and thus

$$|\Delta_j f(x)| \leq |\Delta_j f| * |K_j|(x) \leq C|\Delta_j f| * T_j(x).$$

In view of (4.1), this gives the desired conclusion $|\Delta_j f(x)| \lesssim \omega_j(x)$.     □

We also have:

**Proposition 4.3.** *(1) For all $j \in \mathbb{Z}$,*

$$\omega_j \lesssim 2^{\kappa\sigma} \mathcal{M}\mathcal{M}\Delta_j f. \tag{4.6}$$

*(2) For all $j \in \mathbb{Z}$,*

$$\|\omega_j\|_{L^\infty} \lesssim 2^{\kappa\sigma} \|\Delta_j f\|_{L^\infty} \lesssim 2^{\kappa\sigma} \|f\|_{\dot{F}_q^{\alpha,p}}. \tag{4.7}$$

*Proof.* From (4.1), we deduce

$$\omega_j(x) \lesssim \sum_{r \in \mathbb{Z}^d} T_j|\Delta_j f|(x + 2^{-j}r) E(r).$$

Using (4.4) and (4.5), we get

$$\omega_j(x) \lesssim \sum_{r \in \mathbb{Z}^d} \int_{(0,1)^d} T_j|\Delta_j f|(x + 2^{-j}(r+y)) E(r+y) \, dy = E_j T_j|\Delta_j f|(x).$$

We then observe that $T_j|\Delta_j f|(x) \lesssim \mathcal{M}|\Delta_j f|(x)$ and thus $E_j T_j|\Delta_j f|(x) \lesssim 2^{\kappa\sigma} \mathcal{M}(T_j|\Delta_j f|)$. Both are consequences of (2.9). This proves the first item. The second item is now an easy consequence of (4.6) and the Bernstein inequality (2.12).     □

The derivative estimates for $\omega_j$ can be obtained in a similar manner to the proofs of [33, Proposition 9.6] and [33, Proposition 9.7] and they depend on how many derivatives are computed in the good directions $x' := (x_1, \ldots, x_\kappa)$.

**Proposition 4.4.** *For every $\gamma' \in \mathbb{N}^\kappa$ and $\gamma \in \mathbb{N}^d$,*

$$|\partial^\gamma \partial_{x'}^{\gamma'} \omega_j| \lesssim 2^{(|\gamma'|+|\gamma|)j} 2^{-|\gamma'|\sigma} \omega_j. \tag{4.8}$$

*The implicit constant may depend on $|\gamma|$ and $|\gamma'|$ but neither on $j$ nor on $f$.*

*Proof.* We can assume that $\omega_j \not\equiv 0$, and thus $\omega_j(x) > 0$ for all $x \in \mathbb{R}^d$, see Remark 3.1. The function $u : x \mapsto (1+\|x\|^2)^{1/2}$ is the Euclidean norm of the vector $(1,x)$ in $\mathbb{R}^{d+1}$. By homogeneity, it follows that $|\partial^\gamma u(x)| \lesssim 1$ for all $\gamma \in \mathbb{N}^d \setminus \{0\}$. By the Faà di Bruno formula (i.e. the expression of the iterated derivatives of the composition of two functions), we obtain the pointwise estimate $|\partial^\gamma(\exp \circ (pu))(x)| \lesssim \exp \circ (pu)(x)$. By definition of $E$, see (3.1), it follows that for every $\gamma \in \mathbb{N}^d$, $\gamma' \in \mathbb{N}^\kappa$,

$$\left|\partial^\gamma \partial^{\gamma'}_{x'} E^p(2^j x - r)\right| \lesssim 2^{j(|\gamma|+|\gamma'|)} 2^{-|\gamma'|\sigma} E^p(2^j x - r). \tag{4.9}$$

By definition of $\omega_j$, see (3.3), it follows that

$$\left|\partial^\gamma \partial^{\gamma'}_{x'} \omega_j^p\right| \lesssim 2^{j(|\gamma|+|\gamma'|)} 2^{-|\gamma'|\sigma} \omega_j^p.$$

Writing $\omega_j = \left(\omega_j^p\right)^{1/p}$, the Faà di Bruno formula applied to the functions $\omega_j^p$ and $t \mapsto t^{1/p}$ gives

$$\left|\partial^\gamma \partial^{\gamma'}_{x'} \omega_j\right| \lesssim 2^{j(|\gamma|+|\gamma'|)} 2^{-|\gamma'|\sigma} \omega_j.$$

This proves the proposition. $\square$

From Proposition 4.4, we deduce:

**Proposition 4.5.** *For every $\gamma' \in \mathbb{N}^\kappa$ and $\gamma \in \mathbb{N}^d$,*

$$|\partial^\gamma \partial^{\gamma'}_{x'} \zeta_j| \lesssim 2^{(|\gamma'|+|\gamma|)j} 2^{-|\gamma'|\sigma}.$$

*Proof.* Since the result is obvious when $\zeta_j = 0$, we assume that $\zeta_j \neq 0$. We write

$$\zeta_j(x) = \zeta\left(\frac{2^{\alpha j}\omega_j}{v_j}\right)$$

where $v_j(x) = \sum_{k<j, k\equiv j \pmod R} 2^{\alpha k}\omega_k$. By Proposition 4.4,

$$|\partial^\gamma \partial^{\gamma'}_{x'} v_j| \lesssim \sum_{\substack{k<j \\ k\equiv j \pmod R}} 2^{\alpha k} 2^{(|\gamma|+|\gamma'|)k} 2^{-|\gamma'|\sigma} \omega_k \lesssim 2^{(|\gamma|+|\gamma'|)j} 2^{-|\gamma'|\sigma} v_j. \tag{4.10}$$

We now prove by induction on $|\gamma|+|\gamma'|$ that

$$\left|\partial^\gamma \partial^{\gamma'}_{x'} \frac{1}{v_j}\right| \lesssim 2^{(|\gamma|+|\gamma'|)j} 2^{-|\gamma'|\sigma} \frac{1}{v_j}. \tag{4.11}$$

Since $0 = \partial^\gamma \partial^{\gamma'}_{x'}(v_j \cdot (1/v_j))$, the Leibniz formula implies

$$v_j \partial^\gamma \partial^{\gamma'}_{x'} \frac{1}{v_j} = -\sum_{\substack{\beta \leq \gamma, \beta' \leq \gamma' \\ |\beta|+|\beta'|<|\gamma|+|\gamma'|}} \binom{\gamma}{\beta}\binom{\gamma'}{\beta'} \partial^{\gamma-\beta}\partial^{\gamma'-\beta'}_{x'} v_j \, \partial^\beta \partial^{\beta'}_{x'}\left(\frac{1}{v_j}\right).$$

Using (4.10) and (4.11) for every $|\beta|+|\beta'| < |\gamma|+|\gamma'|$, it then follows that

$$\left|\partial^\gamma \partial^{\gamma'}_{x'} \frac{1}{v_j}\right| \lesssim \frac{1}{v_j} \sum_{\substack{\beta \leq \gamma, \beta' \leq \gamma' \\ |\beta|+|\beta'|<|\gamma|+|\gamma'|}} 2^{(|\gamma-\beta|+|\gamma'-\beta'|)j} 2^{-|\gamma'-\beta'|\sigma} v_j \cdot 2^{(|\beta|+|\beta'|)j} 2^{-|\beta'|\sigma} \frac{1}{v_j} \lesssim 2^{(|\gamma|+|\gamma'|)j} 2^{-|\gamma'|\sigma} \frac{1}{v_j}.$$

By Leibniz formula and Proposition 4.4, this gives

$$\left|\partial^\gamma \partial^{\gamma'}_{x'}\left(\frac{2^{\alpha j}\omega_j}{v_j}\right)\right| \lesssim 2^{(|\gamma|+|\gamma'|)j} 2^{-|\gamma'|\sigma}\left(\frac{2^{\alpha j}\omega_j}{v_j}\right)$$

and the desired estimate now follows from the Faà di Bruno formula applied to the functions $\zeta$ and $\frac{2^{\alpha j}\omega_j}{v_j}$, and also the fact that $\zeta_j(x) = 0$ when $2^{\alpha j}\omega_j(x) > v_j(x)$.

□

4.2. **Integral estimates.** We first establish:

**Proposition 4.6.** *For $1 < p, q < \infty$, $\alpha > 0$,*

(4.12) $$\|\|(2^{\alpha j}\omega_j)\|_{\ell^q(j)}\|_{L^p} \lesssim 2^{\kappa\sigma}\|f\|_{\dot{F}_q^{\alpha,p}}$$

*Proof.* This follows from item 1 in Proposition 4.3 and (2.8). □

The key result of this section is an integral estimate on $\sup_j(2^{\alpha j}\omega_j)$, which will be used crucially to bound $\|\partial_1(h - \tilde{h})\|_{\dot{F}_q^{\alpha,p}}$ in section 5.

**Proposition 4.7.** *One has*

(4.13) $$\|\sup_{j\in\mathbb{Z}}(2^{\alpha j}\omega_j)\|_{L^p} \lesssim \sigma 2^{\frac{\kappa\sigma}{p}}\|f\|_{\dot{F}_q^{\alpha,p}}.$$

The proof of Proposition 4.7 is more involved than the previous ones. It relies on the following estimate:

**Proposition 4.8.** *Let $p \in (1, \infty)$ and $q \in [1, \infty]$. Then there exists $C = C(p, q, d) > 0$ such that for every $f = (f_j)_{j\in\mathbb{Z}} \in L^p(\mathbb{R}^d; \ell^q(\mathbb{Z}))$, for every $r \in \mathbb{R}^d$,*

(4.14) $$\|\|(T_j|f_j|(\cdot + 2^{-j}r))_{j\in\mathbb{Z}}\|_{\ell^q(\mathbb{Z})}\|_{L^p(\mathbb{R}^d)} \leq C\ln(2 + \|r\|)\|\|(f_j)_{j\in\mathbb{Z}}\|_{\ell^q(\mathbb{Z})}\|_{L^p(\mathbb{R}^d)}.$$

The proof of Proposition 4.8 will be given in Appendix 10 below. Let us now derive Proposition 4.7 from Proposition 4.8.

*Proof.* By Proposition 4.1, for every $x \in \mathbb{R}^d$,

$$\sup_{j\in\mathbb{Z}}(2^{\alpha j}\omega_j(x))^p \lesssim \sup_{j\in\mathbb{Z}}\sum_{r\in\mathbb{Z}^d}(2^{\alpha j}T_j|\Delta_j f|(x + 2^{-j}r)E(r))^p$$

$$\lesssim \sum_{r\in\mathbb{Z}^d} E(r)^p \sup_{j\in\mathbb{Z}}(2^{\alpha j}T_j|\Delta_j f|(x + 2^{-j}r))^p$$

$$= \sum_{r\in\mathbb{Z}^d} E(r)^p \left[\sup_{j\in\mathbb{Z}}(2^{\alpha j}T_j|\Delta_j f|(x + 2^{-j}r))\right]^p$$

$$\lesssim \sum_{r\in\mathbb{Z}^d} E(r)^p \|(2^{\alpha j}T_j|\Delta_j f|(x + 2^{-j}r))_{j\in\mathbb{Z}}\|_{\ell^q}^p.$$

Integrating over $x \in \mathbb{R}^d$, we get

$$\|\sup_{j\in\mathbb{Z}} 2^{\alpha j}\omega_j\|_{L^p}^p \lesssim \sum_{r\in\mathbb{Z}^d} E(r)^p\|\|(2^{\alpha j}T_j|\Delta_j f|(\cdot + 2^{-j}r))_{j\in\mathbb{Z}}\|_{\ell^q}\|_{L^p}^p.$$

By Proposition 4.8, this gives

$$\|\sup_{j\in\mathbb{Z}} 2^{\alpha j}\omega_j\|_{L^p}^p \lesssim \sum_{r\in\mathbb{Z}^d} E(r)^p[\ln(2 + \|r\|)]^p\|f\|_{\dot{F}_q^{\alpha,p}}^p.$$

In order to estimate $\sum_{r\in\mathbb{Z}^d} E(r)^p[\ln(2 + \|r\|)]^p$, we first observe that for every $r \in \mathbb{Z}^d$, for every $x \in r + [0,1]^d$,

$$E(r) \leq e^{-\|r_\sigma\|} \leq e^{-\|r_\sigma\|_1/\sqrt{d}} \leq e^{(2^{-\sigma}\kappa + (d-\kappa))/\sqrt{d}}e^{-\|x_\sigma\|_1/\sqrt{d}}.$$

Here $\|x\|_1$ is the $\ell^1$-norm given by $|x_1| + \cdots + |x_d|$. Moreover, $\ln(2 + \|r\|) \leq \ln(3 + \|x\|_1)$. It follows that

$$\sum_{r \in \mathbb{Z}^d} E(r)^p (\ln(2 + \|r\|))^p \lesssim \int_{\mathbb{R}^d} e^{-p\|x_\sigma\|_1/\sqrt{d}} (\ln(3 + \|x\|_1))^p \, dx$$

$$\leq \int_{\mathbb{R}^d} e^{-p\|x_\sigma\|_1/\sqrt{d}} (\ln(3 + 2^\sigma \|x_\sigma\|_1))^p \, dx$$

$$= 2^{\sigma\kappa} \int_{\mathbb{R}^d} e^{-p\|x\|_1/\sqrt{d}} [\ln(3 + 2^\sigma \|x\|_1)]^p \, dx$$

$$\lesssim \sigma^p 2^{\sigma\kappa}.$$

This completes the proof of Proposition 4.7. □

What will be important for us above is that the power of $2^\sigma$ in (4.13), namely $\frac{\kappa}{p}$, is strictly less than 1.

We will use Proposition 4.7 in the following form:

**Lemma 4.9.**

(4.15) $$\left\| \left\| 2^{\alpha m} \omega_m \chi_{2^{\alpha m} \omega_m > \frac{1}{2} \sum_{k<m, k \equiv m \pmod{R}} 2^{\alpha k} \omega_k} \right\|_{\ell^q(m)} \right\|_{L^p} \lesssim R\sigma 2^{\frac{\kappa\sigma}{p}} \|f\|_{\dot{F}_q^{\alpha,p}}.$$

*Proof.* Since $\ell^1(\mathbb{Z})$ continuously embeds in $\ell^q(\mathbb{Z})$, one gets

$$\left\| \left\| 2^{\alpha m} \omega_m \chi_{2^{\alpha m} \omega_m > \frac{1}{2} \sum_{k<m, k \equiv m \pmod{R}} 2^{\alpha k} \omega_k} \right\|_{\ell^q(m)} \right\|_{L^p} \lesssim \left\| \sum_{m \in \mathbb{Z}} 2^{\alpha m} \omega_m \chi_{2^{\alpha m} \omega_m > \frac{1}{2} \sum_{k<m, k \equiv m \pmod{R}} 2^{\alpha k} \omega_k} \right\|_{L^p}.$$

It is enough to combine (3.8) with Proposition 4.7 to conclude the proof of the lemma. □

We end up this section by establishing the expected $L^\infty$ bounds on $F$ under a smallness condition on $\|f\|_{\dot{F}_q^{\alpha,p}}$. More precisely, in view of (4.2) and (4.7), there exists $\eta > 0$ such that if $2^{\kappa\sigma} \|f\|_{\dot{F}_q^{\alpha,p}} \leq \eta$, then for every $j \in \mathbb{Z}$

$$\|\omega_j\|_{L^\infty}, \|\Delta_j f\|_{L^\infty} < 1.$$

By definition of $U_j, h_j$ and $g_j$, see (3.4) and (3.5), this implies

(4.16) $$\|U_j\|_{L^\infty}, \|h_j\|_{L^\infty}, \|g_j\|_{L^\infty} < 1.$$

We can also obtain $L^\infty$ bounds on $G_j, \tilde{h}$ and $\tilde{g}$:

**Lemma 4.10.** *Assume that $2^{\kappa\sigma} \|f\|_{\dot{F}_q^{\alpha,p}} \leq \eta$ with $\eta$ as above. Then*

(1) $\left\|\tilde{h}\right\|_{L^\infty} \lesssim 1$,

(2) *there exists $j_0 \in \mathbb{N}$ such that for every $x \in \mathbb{R}^d$, and every $j \in \mathbb{Z}$,*

$$|G_j(x)| \leq \frac{\min\left(2^{-\alpha R}, 2^{-\alpha(j-j_0)}\right)}{1 - 2^{-\alpha R}}.$$

*In particular, $\|G_j\|_{L^\infty} \leq \frac{1}{2^{\alpha R}-1}$.*

(3) *The infinite products involved in the definition (3.6) of $\tilde{g}$ are uniformly convergent on $\mathbb{R}^d$. If we further assume that $\alpha R > 1$, then $\|G_j\|_{L^\infty} < 1$ and $\|\tilde{g}\|_{L^\infty} \lesssim R$.*

20 PIERRE BOUSQUET, EMMANUEL RUSS, YI WANG, PO-LAM YUNG*Proof.* Using that

$$|h_j(x)| = (1 - \zeta_j(x)) |\Delta_j f(x)| \lesssim (1 - \zeta_j(x))\omega_j(x) = U_j(x),$$

and that $0 \leq U_j \leq 1$ by the choice of $\eta$, we have

$$\left\|\tilde{h}\right\|_{L^\infty} \lesssim \sum_{j=-\infty}^{\infty} U_j \prod_{j'>j} (1 - U_{j'})$$

which implies the first item by (3.7). We now estimate $G_j$. Let $j_0 \in \mathbb{N}$ be an index for which $\Delta_j f \equiv 0$ for all $j > j_0$. Then $\omega_j \equiv 0$ for all $j > j_0$. By the choice of $\eta$, $\|\omega_j\|_{L^\infty} < 1$ for every $j \in \mathbb{Z}$. It follows that for every $x \in \mathbb{R}^d$,

$$0 \leq G_j(x) < \sum_{\substack{t>0, j-t\leq j_0 \\ t\equiv 0 (\mathrm{mod}\ R)}} 2^{-\alpha t} \leq \sum_{k\geq k_0} 2^{-\alpha R k}$$

where $k_0$ is the lowest positive integer such $k_0 R \geq j - j_0$. This implies

$$G_j(x) \leq \frac{2^{-\alpha R k_0}}{1 - 2^{-\alpha R}} \leq \frac{\min\left(2^{-\alpha R}, 2^{-\alpha(j-j_0)}\right)}{1 - 2^{-\alpha R}}$$

and the second item follows.

Moreover, whenever $j > j_0$, $\|G_j\|_{L^\infty(\mathbb{R}^d)} \lesssim 2^{-\alpha(j-j_0)}$ (with an implicit constant depending on $R$) from which we obtain the uniform convergence of $\displaystyle\prod_{\substack{j'>j \\ j'\equiv c (\mathrm{mod}\ R)}} (1 - G_{j'})$ on $\mathbb{R}^d$ for all $j$. This implies the first part of the third item. Finally, in order to obtain the estimate for $\tilde{g}$, we assume that $\alpha R > 1$. By the second item, this implies $\|G_j\|_{L^\infty} < 1$. We next observe that when $\zeta_j(x) > 0$,

$$2^{\alpha j}\omega_j(x) \leq \sum_{\substack{k<j \\ k\equiv j (\mathrm{mod}\ R)}} 2^{\alpha k}\omega_k$$

and thus

(4.17) $$|g_j(x)| \lesssim \zeta_j(x)\omega_j(x) \lesssim \sum_{\substack{k<j \\ k\equiv j (\mathrm{mod}\ R)}} 2^{\alpha(k-j)}\omega_k = G_j.$$

It then follows that

$$\|\tilde{g}\|_{L^\infty} \lesssim \sum_{c=0}^{R-1} \sum_{j\equiv c (\mathrm{mod}\ R)} G_j \prod_{\substack{j'>j \\ j'\equiv c (\mathrm{mod}\ R)}} (1 - G_{j'}) \lesssim R.$$

This completes the proof. □

*In the next two sections, we will always assume that*

(4.18) $$\|f\|_{\dot{F}_q^{\alpha,p}} \leq 2^{-\kappa\sigma}\eta$$

and also that

(4.19) $$\alpha R > 1.$$

## 5. Estimating $h - \tilde{h}$

We still write $\partial_{x'}$ for a derivative in any of the "good" directions, namely $\partial_1, \ldots, \partial_\kappa$. This section is devoted to the proof of the $\dot{F}_q^{\alpha-1,p}$ estimate for the derivatives of $h - \tilde{h}$:

**Proposition 5.1.** *Let $\alpha, p, q$ and $\kappa$ be as in Theorem 1.1. Define $h$ (resp. $\tilde{h}$) by (3.4) (resp. (3.5)). Then*

$$\|\partial_{x'}(h-\tilde{h})\|_{\dot{F}_q^{\alpha-1,p}} \lesssim R\sigma^2 \left( 2^{(-\min(1,\alpha)+\frac{\kappa}{p})\sigma}\|f\|_{\dot{F}_q^{\alpha,p}} + 2^{\left(\max(1-\alpha,0)+\kappa\left(1+[\alpha]+\frac{1}{p}\right)\right)\sigma}\|f\|_{\dot{F}_q^{\alpha,p}}^2 \right),$$

*and for any $1 \leq i \leq d$,*

$$\|\partial_i(h-\tilde{h})\|_{\dot{F}_q^{\alpha-1,p}} \lesssim R\sigma^2 \left( 2^{\frac{\kappa}{p}\sigma}\|f\|_{\dot{F}_q^{\alpha,p}} + 2^{\left(\max(1-\alpha,0)+\kappa\left(1+[\alpha]+\frac{1}{p}\right)\right)\sigma}\|f\|_{\dot{F}_q^{\alpha,p}}^2 \right),$$

*where the implicit constants depend on $\alpha$, $p$ and $q$ but not on $f, R, \sigma$.*

Here, $[\alpha]$ is the integer part of $\alpha$.

*Proof of Proposition 5.1:* We only prove the first inequality of the statement. The proof of the second one is very similar (and easier to establish).

**Step 1.** Let

$$V_j := \sum_{j'<j} h_{j'} \prod_{j'<j''<j} (1 - U_{j''}).$$

Then

(5.1) $$h - \tilde{h} = \sum_j U_j V_j.$$

Identity (5.1) is a consequence of Lemma 3.2:

$$h - \tilde{h} = \sum_{j=-\infty}^{\infty} h_j \left( 1 - \prod_{j'>j}(1-U_{j'}) \right)$$
$$= \sum_j h_j \sum_{j'>j} U_{j'} \prod_{j<j''<j'} (1-U_{j''})$$
$$= \sum_{j'} U_{j'} \sum_{j<j'} h_j \prod_{j<j''<j'} (1-U_{j''}) = \sum_j U_j V_j.$$

**Step 2: estimates on $U_j$, $V_j$ and their derivatives.** Let us first collect the estimates for $U_j$:

**Lemma 5.2.**   (1) *For every $\gamma' \in \mathbb{N}^\kappa, \gamma \in \mathbb{N}^d$ and every $m \in \mathbb{Z}, x \in \mathbb{R}^d$,*

$$|\partial^\gamma \partial_{x'}^{\gamma'} U_m(x)| \lesssim 2^{m|\gamma|} 2^{(m-\sigma)|\gamma'|} \omega_m(x) \chi_{2^{\alpha m}\omega_m > \frac{1}{2} \sum_{k<m, k\equiv m(\mathrm{mod}\ R)} 2^{\alpha k}\omega_k}(x),$$

(2) *For every $\gamma \in \mathbb{N}^d$, $\|\partial^\gamma U_m\|_{L^\infty} \lesssim 2^{m|\gamma|}2^{\kappa\sigma}\|f\|_{\dot{F}_q^{\alpha,p}}$.*

In the above statement as well as in the lemmata below, the implicit constants may depend on the number of derivatives $|\gamma|$ and $|\gamma'|$ (but neither on $m, x$ nor on $f$).

*Proof.* When there exists $k < m$ with $k \equiv m \pmod{R}$ such that $\omega_k \not\equiv 0$, estimate (1) follows from the definition of the functions $U_m$, see (3.5), Proposition 4.4, Proposition 4.5 and the Leibniz rule. We also rely on the fact that $(1-\zeta_m) \equiv 0$ and thus $U_m \equiv 0$ on the set where $2^{\alpha m}\omega_m \leq \frac{1}{2}\sum_{k<m, k\equiv m \pmod{R}} 2^{\alpha k}\omega_k$. When $\omega_k \equiv 0$ for every $k < m, k \equiv m \pmod{R}$, $U_m = \omega_m$ and one therefore has
$$|\partial^\gamma \partial_{x'}^{\gamma'} U_m(x)| \lesssim 2^{m|\gamma|} 2^{(m-\sigma)|\gamma'|} \omega_m(x).$$
If $\omega_m(x) > 0$, the conclusion readily follows. Otherwise, $\omega_m$ identically vanishes and the estimate is obvious.

It follows from the first item and (4.7) that, for every $\gamma \in \mathbb{N}^d$,
$$\|\partial^\gamma U_m\|_{L^\infty} \lesssim 2^{m|\gamma|}\|\omega_m\|_{L^\infty} \lesssim 2^{m|\gamma|} 2^{\kappa\sigma}\|f\|_{\dot{F}_q^{\alpha,p}},$$
which proves the second item. $\square$

**Lemma 5.3.** *For all $m \in \mathbb{Z}$ and $\gamma \in \mathbb{N}^d$,*
$$\|\partial^\gamma h_m\|_{L^\infty} \lesssim 2^{m|\gamma|}\|f\|_{\dot{F}_q^{\alpha,p}}.$$

*Proof.* By definition of $h_m$, the Leibniz rule and Proposition 4.5,
$$\|\partial^\gamma h_m\|_{L^\infty} \lesssim \sum_{0 \leq \gamma' \leq \gamma} 2^{m|\gamma-\gamma'|} \left\|\partial^{\gamma'}\Delta_m f\right\|_{L^\infty}.$$
We now rely on (2.11) to get
$$\|\partial^\gamma h_m\|_{L^\infty} \lesssim \sum_{0 \leq \gamma' \leq \gamma} 2^{m|\gamma-\gamma'|} 2^{m|\gamma'|}\|f\|_{\dot{F}_q^{\alpha,p}} \lesssim 2^{m|\gamma|}\|f\|_{\dot{F}_q^{\alpha,p}}.$$
$\square$

Here are now the estimates for $V_m$:

**Lemma 5.4.** (1) *For every $m \in \mathbb{Z}$, $|V_m| \lesssim 1$.*
(2) *For every $\gamma \in \mathbb{N}^d$, $|\partial^\gamma V_m| \lesssim 2^{m|\gamma|} 2^{\kappa|\gamma|\sigma}\|f\|_{\dot{F}_q^{\alpha,p}}$.*

*Proof.* The first item follows from the construction of $V_m$, (3.7) and the fact that for all $x \in \mathbb{R}^d$,
$$|h_j(x)| \lesssim (1-\zeta_j(x))\omega_j(x) = U_j(x) \leq 1,$$
where the last inequality above is due to (4.16).

Let us prove the second item. Arguing as in [33, (6.8)], one obtains
$$\partial_k V_m = \sum_{m'<m} (\partial_k h_{m'} - V_{m'}\partial_k U_{m'}) \prod_{m'<m''<m} (1-U_{m''}).$$
Using this calculation, one can prove by induction on $|\gamma|$, $\gamma \in \mathbb{N}^d$, that
$$(5.2) \qquad \partial^\gamma V_m = \sum_{m'<m} \left(\partial^\gamma h_{m'} - \sum_{0<\alpha\leq\gamma} c_{\alpha,\gamma}\partial^\alpha U_{m'}\partial^{\gamma-\alpha}V_{m'}\right) \prod_{m'<m''<m} (1-U_{m''})$$

where $c_{\alpha,\gamma}$ is some positive integer for each $0 < \alpha < \gamma$. Indeed, for any finite sequence $I = (I_1, \ldots, I_k)$ where $I_1, \ldots, I_k \in \{1, \ldots, d\}$ and $k \in \mathbb{N}$, we have

$$\partial_I V_m = \sum_{m'<m} \left( \partial_I h_{m'} - \sum_{\substack{J \neq \emptyset \\ J \text{ subsequence of } I}} \partial_J U_{m'} \, \partial_{I \setminus J} V_{m'} \right) \prod_{m'<m''<m} (1 - U_{m''})$$

where for each non-empty subsequence $J$ of $I$, $I \setminus J$ is the subsequence of $I$ obtained by removing $J$ from $I$. (A subsequence $J$ of $I = (I_1, \ldots, I_k)$ is a finite sequence of the form $(I_{i_1}, \ldots, I_{i_\ell})$ with $\ell \leq k$ and $1 \leq i_1 < \cdots < i_\ell \leq k$.) From (5.2) we deduce that for every $\gamma \in \mathbb{N}^d$,

$$\|\partial^\gamma V_m\|_{L^\infty} \leq \sum_{m'<m} \left( \|\partial^\gamma h_{m'}\|_{L^\infty} + C \sum_{0 < \alpha \leq \gamma} \|\partial^\alpha U_{m'}\|_{L^\infty} \|\partial^{\gamma-\alpha} V_{m'}\|_{L^\infty} \right)$$

By item (2) in Lemma 5.2,
$$\|\partial^\alpha U_{m'}\|_{L^\infty} \lesssim 2^{m'|\alpha|} 2^{\kappa\sigma} \|f\|_{\dot{F}_q^{\alpha,p}}.$$

Moreover, as a consequence of Lemma 5.3, $\|\partial^\gamma h_{m'}\|_{L^\infty} \lesssim 2^{m'|\gamma|} \|f\|_{\dot{F}_q^{\alpha,p}}$. This implies

$$\|\partial^\gamma V_m\|_{L^\infty} \lesssim \sum_{m'<m} \left( 2^{m'|\gamma|} + \sum_{0<\alpha\leq\gamma} 2^{m'|\alpha|} 2^{\kappa\sigma} \|\partial^{\gamma-\alpha} V_{m'}\|_{L^\infty} \right) \|f\|_{\dot{F}_q^{\alpha,p}}.$$

The result then follows by induction on $|\gamma|$, since induction hypothesis implies $\|\partial^{\gamma-\alpha} V_{m'}\|_{L^\infty} \lesssim 2^{m'(|\gamma|-|\alpha|)} 2^{(|\gamma|-|\alpha|)\kappa\sigma}$ if $0 < \alpha \leq \gamma$. □

**Step 3 Completion of the proof of Proposition 5.1.** From (5.1) we see that

$$\begin{aligned}
\|\partial_{x'}(h - \tilde{h})\|_{\dot{F}_q^{\alpha-1,p}} &= \|\|2^{(\alpha-1)m} \Delta_m (\partial_{x'}(h - \tilde{h}))\|_{\ell^q(m)}\|_{L^p} \\
&= \|\|2^{(\alpha-1)m} \Delta_m (\partial_{x'} \sum_{j \in \mathbb{Z}} U_j V_j)\|_{\ell^q(m)}\|_{L^p} \\
&= \|\|2^{(\alpha-1)m} \Delta_m (\partial_{x'} \sum_{r \in \mathbb{Z}} U_{r+m} V_{r+m})\|_{\ell^q(m)}\|_{L^p}.
\end{aligned}$$

By the triangle inequality, this gives

(5.3) $$\|\partial_{x'}(h - \tilde{h})\|_{\dot{F}_q^{\alpha-1,p}} \leq \sum_{r \in \mathbb{Z}} \|\|2^{(\alpha-1)m} \Delta_m (\partial_{x'}(U_{r+m} V_{r+m}))\|_{\ell^q(m)}\|_{L^p}.$$

We now introduce a positive integer $T$ to be defined later and we split the sum into three parts that we estimate separately: $\sum_{r>T} \ldots$, $\sum_{0 \leq r \leq T} \ldots$ and $\sum_{r<0} \ldots$.

**5.1. Estimate of $\sum_{r>T}$.** In this case, we let $\partial_{x'}$ differentiate the Littlewood-Paley projection $\Delta_m$. By Lemma 9.1 below and the fact that $\|V_m\|_{L^\infty} \lesssim 1$ (Lemma 5.4),

$$\|\|2^{(\alpha-1)m} \Delta_m (\partial_{x'}(U_{r+m} V_{r+m}))\|_{\ell^q(m)}\|_{L^p} \lesssim \|\|2^{(\alpha-1)m} 2^m U_{r+m} V_{r+m}\|_{\ell^q(m)}\|_{L^p}$$

(5.4) $$\lesssim 2^{-\alpha r} \|\|2^{\alpha m} U_m\|_{\ell^q(m)}\|_{L^p}.$$

By Lemma 4.9, and the definition of $U_j$,

$$\||2^{\alpha m}U_m\|_{\ell^q(m)}\|_{L^p} \lesssim \||2^{\alpha m}\omega_m \chi_{2^{\alpha m}\omega_m > \frac{1}{2}\sum_{k<m, k\equiv m\,(\text{mod }R)}2^{\alpha k}\omega_k}\|_{\ell^q(m)}\|_{L^p}$$

$$\lesssim R\sigma 2^{\frac{\kappa\sigma}{p}}\|f\|_{\dot{F}_q^{\alpha,p}}.$$

Hence, by summing over $r > T$, one gets

$$\sum_{r>T}\||2^{(\alpha-1)m}\Delta_m(\partial_{x'}(U_{r+m}V_{r+m}))\|_{\ell^q(m)}\|_{L^p} \lesssim R\sigma 2^{-\alpha T + \frac{\kappa\sigma}{p}}\|f\|_{\dot{F}_q^{\alpha,p}}.$$

**5.2. Estimate of $\sum_{r<0}$.** Let $a$ be the integer part of $\alpha$. By Lemma 9.2 and (2.3), there exist Schwartz functions $\Delta^{(\gamma)}$ such that

$$\Delta = \sum_{|\gamma|=a}\partial^\gamma \Delta^{(\gamma)}.$$

Then

$$\Delta_m(x) = 2^{md}\Delta(2^m x) = \sum_{|\gamma|=a}2^{md}[\partial^\gamma \Delta^{(\gamma)}](2^m x) = 2^{-ma}\sum_{|\gamma|=a}\partial^\gamma[(\Delta^{(\gamma)})_m](x)$$

where $(\Delta^{(\gamma)})_m(x) = 2^{md}\Delta^{(\gamma)}(2^m x)$. Hence,

$$\||2^{(\alpha-1)m}\Delta_m(\partial_{x'}(U_{m+r}V_{m+r}))\|_{\ell^q(m)}\|_{L^p}$$

(5.5)
$$\leq \sum_{|\gamma|=a}\||2^{(\alpha-1-a)m}\partial^\gamma[(\Delta^{(\gamma)})_m](\partial_{x'}(U_{m+r}V_{m+r}))\|_{\ell^q(m)}\|_{L^p}$$

$$= \sum_{|\gamma|=a}\||2^{(\alpha-1-a)m}(\Delta^{(\gamma)})_m(\partial^{(\gamma)}\partial_{x'}(U_{m+r}V_{m+r}))\|_{\ell^q(m)}\|_{L^p}$$

$$\lesssim \sum_{|\gamma|=a}\||2^{(\alpha-1-a)m}(\partial^\gamma \partial_{x'}(U_{m+r}V_{m+r}))\|_{\ell^q(m)}\|_{L^p}.$$

In the last line, we have used Lemma 9.1, applied to the function $\Delta^{(\gamma)}$. This implies
(5.6)
$$\||2^{(\alpha-1)m}\Delta_m(\partial_{x'}(U_{m+r}V_{m+r}))\|_{\ell^q(m)}\|_{L^p} \lesssim 2^{-(\alpha-1-a)r}\sum_{|\gamma|=a}\||2^{(\alpha-1-a)m}(\partial^\gamma \partial_{x'}(U_m V_m))\|_{\ell^q(m)}\|_{L^p}.$$

Now, by the Leibniz rule,

(5.7)
$$|\partial^\gamma \partial_{x'}(U_m V_m)| \lesssim |V_m(\partial^\gamma \partial_{x'}U_m)| + \sum_{\ell=0}^{a}|(\partial_x^\ell U_m)(\partial_x^{a+1-\ell}V_m)|.$$

Here, $\partial_x^\ell$ refers to the full partial differential operator of order $\ell$ and similarly for $\partial_x^{a+1-\ell}$. By Lemmata 5.2 and 5.4, one gets

$$|\partial^\gamma \partial_{x'}(U_m V_m)| \lesssim 2^{-\sigma}2^{m(a+1)}\omega_m \chi_{2^{\alpha m}\omega_m > \frac{1}{2}\sum_{k<m, k\equiv m\,(\text{mod }R)}2^{\alpha k}\omega_k}$$

$$+ \sum_{\ell=0}^{a}2^{m\ell}\omega_m \chi_{2^{\alpha m}\omega_m > \frac{1}{2}\sum_{k<m, k\equiv m\,(\text{mod }R)}2^{\alpha k}\omega_k}2^{(a+1-\ell)m}2^{(a+1-\ell)\kappa\sigma}\|f\|_{\dot{F}_q^{\alpha,p}}$$

$$\lesssim (2^{-\sigma} + 2^{(a+1)\kappa\sigma}\|f\|_{\dot{F}_q^{\alpha,p}})2^{m(a+1)}\omega_m \chi_{2^{\alpha m}\omega_m > \frac{1}{2}\sum_{k<m, k\equiv m\,(\text{mod }R)}2^{\alpha k}\omega_k}.$$

We deduce from (5.6) that

$$\|\|2^{(\alpha-1)m}\Delta_m(\partial_{x'}(U_{m+r}V_{m+r}))\|_{\ell^q(m)}\|_{L^p}$$
$$\lesssim (2^{-\sigma} + 2^{(a+1)\kappa\sigma}\|f\|_{\dot{F}_q^{\alpha,p}})2^{-(\alpha-1-a)r}\|\|2^{\alpha m}\omega_m\chi_{2^{\alpha m}\omega_m > \frac{1}{2}}\sum_{k<m, k\equiv m \pmod{R}} 2^{\alpha k}\omega_k\|_{\ell^q(m)}\|_{L^p}$$
$$\lesssim R\sigma(2^{-\sigma} + 2^{(a+1)\kappa\sigma}\|f\|_{\dot{F}_q^{\alpha,p}})2^{-(\alpha-1-a)r}2^{\frac{\kappa\sigma}{p}}\|f\|_{\dot{F}_q^{\alpha,p}},$$

where the last line follows from Lemma 4.9.

Thus by summing over $r < 0$ (taking into account that $\alpha < 1 + a$),
(5.8)
$$\sum_{r<0} \|\|2^{(\alpha-1)m}\Delta_m(\partial_{x'}(U_{r+m}V_{r+m}))\|_{\ell^q(m)}\|_{L^p} \lesssim R\sigma\left(2^{(-1+\frac{\kappa}{p})\sigma}\|f\|_{\dot{F}_q^{\alpha,p}} + 2^{\kappa\sigma\left(a+1+\frac{1}{p}\right)}\|f\|_{\dot{F}_q^{\alpha,p}}^2\right).$$

5.3. **Estimate of $\sum_{0 \leq r \leq T}$.** This is exactly the same calculation as in the case $\sum_{r \leq 0}$ except that in (5.5) we take $a = 0$; that is, we do not perform the preliminary integration by parts and we keep $\Delta_m$ instead of introducing $\Delta_m^{(\gamma)}$. Hence, when summing over $0 \leq r \leq T$, (5.8) is replaced by
(5.9)
$$\sum_{0 \leq r \leq T} \|\|2^{(\alpha-1)m}\Delta_m(\partial_{x'}(U_{r+m}V_{r+m}))\|_{\ell^q(m)}\|_{L^p} \leq C_\alpha(T)R\sigma\left(2^{(-1+\frac{\kappa}{p})\sigma}\|f\|_{\dot{F}_q^{\alpha,p}} + 2^{\kappa\sigma\left(1+\frac{1}{p}\right)}\|f\|_{\dot{F}_q^{\alpha,p}}^2\right),$$

where

$$C_\alpha(T) = \begin{cases} C2^{(1-\alpha)T} & \text{if } \alpha < 1, \\ CT & \text{if } \alpha = 1, \\ C & \text{if } \alpha > 1. \end{cases}$$

**Remark 5.5.** Note that it is crucial for the sequel to obtain an arbitrarily small factor in front of $\|f\|_{\dot{F}_q^{\alpha,p}}$ in the right-hand sides of (5.8) and (5.9). This in turn follows from Lemma 5.2 and the fact that we take one derivative in a "good" direction in (5.7).

5.4. **Conclusion.** From the three above subsections, one gets, with $a = [\alpha]$,

(1) When $\alpha > 1$, one can take $T = \infty$:
$$\|\partial_{x'}(h - \tilde{h})\|_{\dot{F}_q^{\alpha-1,p}} \lesssim R\sigma\left(2^{(-1+\frac{\kappa}{p})\sigma}\|f\|_{\dot{F}_q^{\alpha,p}} + 2^{\kappa\sigma\left(a+1+\frac{1}{p}\right)}\|f\|_{\dot{F}_q^{\alpha,p}}^2\right).$$

(2) When $\alpha = 1$, one can take $T = \sigma$, which implies
$$\|\partial_{x'}(h - \tilde{h})\|_{\dot{F}_q^{\alpha-1,p}} \lesssim R\sigma\left(\sigma 2^{(-1+\frac{\kappa}{p})\sigma}\|f\|_{\dot{F}_q^{\alpha,p}} + \sigma 2^{\kappa\sigma\left(2+\frac{1}{p}\right)}\|f\|_{\dot{F}_q^{\alpha,p}}^2\right).$$

(3) When $0 < \alpha < 1$, one can take $T = \sigma$, which yields
$$\|\partial_{x'}(h - \tilde{h})\|_{\dot{F}_q^{\alpha-1,p}} \lesssim R\sigma\left(2^{-(\alpha-\frac{\kappa}{p})\sigma}\|f\|_{\dot{F}_q^{\alpha,p}} + 2^{\left(1-\alpha+\kappa\left(1+\frac{1}{p}\right)\right)\sigma}\|f\|_{\dot{F}_q^{\alpha,p}}^2\right).$$

Altogether this proves Proposition 5.1. □

6. Estimating $g - \tilde{g}$

**Proposition 6.1.** *Let $\alpha, p, q$ be as in Theorem 1.1. Define $g$ (resp. $\tilde{g}$) by (3.4) (resp. (3.6)). We also introduce a number $a_\alpha \in (0, \alpha]$ such that $a_\alpha < \frac{\alpha}{1-\alpha}$ when $\alpha < 1$ and $a_\alpha = 1$ when $\alpha \geq 1$. Then*

$$(6.1) \quad \|\partial_x(g - \tilde{g})\|_{\dot{F}_q^{\alpha-1,p}} \lesssim \left(2^{\kappa\sigma} R 2^{-\min(1,\alpha a_\alpha)R} \|f\|_{\dot{F}_q^{\alpha,p}} + 2^{([\alpha]+2)\kappa\sigma} R^2 2^{-\min(1,\alpha a_\alpha)R} \|f\|_{\dot{F}_q^{\alpha,p}}^2\right).$$

*where the implicit constant depends on $\alpha$ and $a_\alpha$, $p$, $q$ but not on $f, R, \sigma$.*

Recall that $[\alpha]$ is the integer part of $\alpha$. Also, we have written $\partial_x$ for any partial derivative $\partial_i$, $1 \leq i \leq d$.

**Remark 6.2.** Note that, contrary to Proposition 5.1, we do not state an improved estimate for the derivatives in the good directions.

**Step 1.** Let

$$H_j := \sum_{\substack{j' < j \\ j' \equiv j (\text{mod } R)}} g_{j'} \prod_{\substack{j' < j'' < j \\ j'' \equiv j (\text{mod } R)}} (1 - G_{j''}).$$

Then, as in Step 1 of the proof of Proposition 5.1,

$$g - \tilde{g} = \sum_j G_j H_j.$$

**Step 2 Estimates on $G_m$, $H_m$ and their derivatives.** Let us collect the upper bounds for $G_m$ and $H_m$ which will be needed in the sequel (see [33, Propositions 11.2-11.7]). In the lemmata below, the implicit constants may depend on $|\gamma|$ but neither on $m$, $\sigma$, $R$ nor on $f$.

**Lemma 6.3.** *For all $m \in \mathbb{Z}$ and all $\gamma \in \mathbb{N}^d$,*

$$|\partial^\gamma G_m| \lesssim 2^{\kappa\sigma} \sum_{\substack{t > 0 \\ t \equiv 0 (\text{mod } R)}} 2^{-\alpha t} 2^{|\gamma|(m-t)} \mathcal{M}\mathcal{M}\Delta_{m-t} f.$$

*Proof.* By definition of $G_m$, see (3.6), and Proposition 4.4, one has

$$|\partial^\gamma G_m| \leq \sum_{\substack{t > 0 \\ t \equiv 0 (\text{mod } R)}} 2^{-\alpha t} |\partial^\gamma \omega_{m-t}| \lesssim \sum_{\substack{t > 0 \\ t \equiv 0 (\text{mod } R)}} 2^{-\alpha t} 2^{|\gamma|(m-t)} \omega_{m-t}.$$

By (4.6), this implies

$$|\partial^\gamma G_m| \lesssim 2^{\kappa\sigma} \sum_{\substack{t > 0 \\ t \equiv 0 (\text{mod } R)}} 2^{-\alpha t} 2^{|\gamma|(m-t)} \mathcal{M}\mathcal{M}\Delta_{m-t} f.$$

□

**Lemma 6.4.** *For all $m \in \mathbb{Z}$ and all $\gamma \in \mathbb{N}^d$,*

$$|\partial^\gamma g_m| \lesssim 2^{|\gamma|m} \mathcal{M}\Delta_m f.$$

*Proof.* By definition of $g_m$ see (3.4), the Leibniz rule and Proposition 4.5,

$$|\partial^\gamma g_m| \lesssim \sum_{0 \leq \gamma' \leq \gamma} 2^{|\gamma-\gamma'|m} \left|\partial^{\gamma'} \Delta_m f\right|.$$

We now rely on (2.10) to get

$$|\partial^\gamma g_m| \lesssim \sum_{0 \leq \gamma' \leq \gamma} 2^{|\gamma-\gamma'|m} 2^{|\gamma'|m} \mathcal{M}\Delta_m f \lesssim 2^{|\gamma|m} \mathcal{M}\Delta_m f.$$

$\square$

**Lemma 6.5.** *For all $m \in \mathbb{Z}$ and all $\gamma \in \mathbb{N}^d$,*

$$|H_m| \lesssim 1 \quad , \quad |\partial^\gamma H_m| \lesssim 2^{|\gamma|\kappa\sigma} \sum_{\substack{t>0 \\ t \equiv 0 (\text{mod } R)}} 2^{|\gamma|(m-t)} \mathcal{M}\mathcal{M}\Delta_{m-t} f.$$

*Proof.* In order to prove the estimate on $|H_m|$, we recall that under the conditions (4.18) and (4.19), $0 \leq G_j \leq 1$. In view of (4.17), namely $|g_j(x)| \lesssim G_j(x)$, this, in conjunction with (3.7), implies $|H_m(x)| \lesssim 1$.

For the second estimate, we argue by induction, assuming the correct bound for $\left|\partial^{\gamma'} H_m\right|$ with $|\gamma'| < |\gamma|$. We have the analogue of (5.2) for $H_m$ in place of $V_m$:

$$\partial^\gamma H_m = \sum_{\substack{m'<m \\ m' \equiv m (\text{mod } R)}} \left(\partial^\gamma g_{m'} - \sum_{0<\gamma'\leq\gamma} c_{\gamma',\gamma} \partial^{\gamma'} G_{m'} \partial^{\gamma-\gamma'} H_{m'}\right) \prod_{\substack{m'<m''<m \\ m'' \equiv m (\text{mod } R)}} (1-G_{m''}).$$

Thus

$$|\partial^\gamma H_m| \lesssim \sum_{\substack{m'<m, \\ m' \equiv m (\text{mod } R)}} |\partial^\gamma g_{m'}| + \sum_{\substack{m'<m, \\ m' \equiv m (\text{mod } R)}} \sum_{0<\gamma'\leq\gamma} \left|\partial^{\gamma'} G_{m'} \partial^{\gamma-\gamma'} H_{m'}\right| = I + II.$$

In view of Lemma 6.4,

$$I \lesssim \sum_{\substack{t>0 \\ t \equiv 0 (\text{mod } R)}} 2^{|\gamma|(m-t)} \mathcal{M}\mathcal{M}(\Delta_{m-t} f).$$

Using instead Lemma 6.3 and the induction assumption, we get

$$II \lesssim 2^{|\gamma|\kappa\sigma} \sum_{\substack{m'<m, \\ m' \equiv m (\text{mod } R)}} \sum_{0<\gamma'\leq\gamma} \sum_{\substack{t>0, l>0, \\ t \equiv l \equiv 0 (\text{mod } R)}} 2^{-\alpha t} 2^{|\gamma'|(m'-t)} \mathcal{M}\mathcal{M}\Delta_{m'-t} f \cdot 2^{|\gamma-\gamma'|(m'-l)} \mathcal{M}\mathcal{M}\Delta_{m'-l} f$$

$$= 2^{|\gamma|\kappa\sigma} \sum_{\substack{m'<m, \\ m' \equiv m (\text{mod } R)}} 2^{|\gamma|m'} \sum_{0<\gamma'\leq\gamma} \sum_{\substack{t>0, l>0, \\ t \equiv l \equiv 0 (\text{mod } R)}} 2^{-\alpha t} 2^{-|\gamma'|t} \mathcal{M}\mathcal{M}\Delta_{m'-t} f \cdot 2^{-|\gamma-\gamma'|l} \mathcal{M}\mathcal{M}\Delta_{m'-l} f.$$

We split the innermost sum into two parts according to whether $t > l$ or $l \geq t$ and estimate by 1 the factor $\mathcal{M}\mathcal{M}\Delta_{m'-t} f$ and $\mathcal{M}\mathcal{M}\Delta_{m'-l} f$ respectively (remember that $\|\Delta_m f\|_{L^\infty} < 1$ under the condition (4.18)); this gives

$$\sum_{\substack{t>0, l>0, \\ t \equiv l \equiv 0 (\text{mod } R)}} 2^{-\alpha t} 2^{-|\gamma'|t} \mathcal{M}\mathcal{M}\Delta_{m'-t} f \cdot 2^{-|\gamma-\gamma'|l} \mathcal{M}\mathcal{M}\Delta_{m'-l} f \lesssim \sum_{\substack{t>0, \\ t \equiv 0 (\text{mod } R)}} 2^{-\alpha t} 2^{-|\gamma|t} \mathcal{M}\mathcal{M}\Delta_{m'-t} f.$$

Finally,

$$II \lesssim 2^{|\gamma|\kappa\sigma} \sum_{\substack{m'<m,\\ m'\equiv m(\text{mod } R)}} 2^{|\gamma|m'} \sum_{\substack{t>0,\\ t\equiv 0(\text{mod } R)}} 2^{-\alpha t} 2^{-|\gamma|t} \mathcal{M}\mathcal{M}\Delta_{m'-t}f$$

$$= 2^{|\gamma|\kappa\sigma} \sum_{\substack{l>0,\\ l\equiv 0(\text{mod } R)}} \sum_{\substack{t>0,\\ t\equiv 0(\text{mod } R)}} 2^{-\alpha t} 2^{|\gamma|(m-l-t)} \mathcal{M}\mathcal{M}\Delta_{m-l-t}f.$$

We sum this double sum by first summing over the pairs $(t,l)$ where $l+t$ is constant, and then sum over the remaining variable. This gives the desired conclusion. $\square$

**Step 3 Completion of the proof.** As in the proof of (5.3),

$$\|\partial_x(g-\tilde{g})\|_{\dot{F}_q^{\alpha-1,p}} \leq \sum_{r\in\mathbb{Z}} \|\|2^{(\alpha-1)m}\Delta_m(\partial_x(G_{r+m}H_{r+m}))\|_{\ell^q(m)}\|_{L^p}.$$

By definition of $G_{r+m}$ (see (3.6)) and the triangle inequality,

$$\|\partial_x(g-\tilde{g})\|_{\dot{F}_q^{\alpha-1,p}} \leq \sum_{\substack{t>0,\\ t\equiv 0(\text{mod } R)}} 2^{-\alpha t} \sum_{r\in\mathbb{Z}} \|\|2^{(\alpha-1)m}\Delta_m(\partial_x(\omega_{r+m-t}H_{r+m}))\|_{\ell^q(m)}\|_{L^p}.$$

For each fixed $t$, we now split the sum over $r$ into three parts as follows

$$\sum_{r\in\mathbb{Z}} \cdots = \sum_{r\leq 0} \cdots + \sum_{0<r<a_\alpha t} \cdots + \sum_{r\geq a_\alpha t} \cdots,$$

where $a_\alpha > 0$ was defined in Proposition 6.1. We proceed to estimate the three terms separately.

6.1. **Estimate of** $\sum_{r\geq a_\alpha t}$. As in the proof of (5.4) we integrate by parts to let $\partial_x$ hit $\Delta_m$, and use $\|H_m\|_{L^\infty} \lesssim 1$ to get

$$\|\|2^{(\alpha-1)m}\Delta_m(\partial_x(\omega_{r+m-t}H_{r+m}))\|_{\ell^q(m)}\|_{L^p} \lesssim 2^{-\alpha(r-t)}\|\|2^{\alpha m}\omega_m\|_{\ell^q(m)}\|_{L^p}.$$

This implies

$$\|\|2^{(\alpha-1)m}\Delta_m(\partial_x(\omega_{r+m-t}H_{r+m}))\|_{\ell^q(m)}\|_{L^p} \lesssim 2^{\kappa\sigma}2^{-\alpha(r-t)}\|f\|_{\dot{F}_q^{\alpha,p}},$$

by Proposition 4.6. Summing over $r$ and $t$, one gets

$$\sum_{\substack{t>0\\ t\equiv 0(\text{mod } R)}} \sum_{r\geq a_\alpha t} 2^{-\alpha t}\|\|2^{(\alpha-1)m}\Delta_m(\partial_x(\omega_{r+m-t}H_{r+m}))\|_{\ell^q(m)}\|_{L^p} \lesssim 2^{\kappa\sigma}\left(\sum_{t\geq R}\sum_{r\geq a_\alpha t} 2^{-\alpha r}\right)\|f\|_{\dot{F}_q^{\alpha,p}}$$

(6.2)
$$\lesssim 2^{\kappa\sigma}2^{-\alpha a_\alpha R}\|f\|_{\dot{F}_q^{\alpha,p}}.$$

6.2. **Estimate of** $\sum_{r\leq 0}$. Let $a \geq 0$ be an integer. Arguing as for the proof of (5.6), we write $\Delta_m$ as an $a$-th derivative and integrate by parts to hit $\partial_x(\omega_{m+r-t}H_{m+r})$, and obtain

(6.3)
$$\|\|2^{(\alpha-1)m}\Delta_m(\partial_x(\omega_{m+r-t}H_{m+r}))\|_{\ell^q(m)}\|_{L^p} \lesssim 2^{-(\alpha-1-a)r}\|\|2^{(\alpha-1-a)m}(\partial_x^{a+1}(\omega_{m-t}H_m))\|_{\ell^q(m)}\|_{L^p}.$$

Here $\partial_x^{a+1}$ refer to the full partial differential operator of order $a+1$.

We now use the fact that for every $\ell \in \mathbb{N}$

$$\omega_{m-t} \lesssim 2^{\kappa\sigma}\mathcal{M}\mathcal{M}(\Delta_{m-t}f) \quad , \quad |\partial_x^\ell \omega_{m-t}| \lesssim 2^{(m-t)\ell}\omega_{m-t}$$

(see Propositions 4.3 and 4.4) and also for $k \in \mathbb{N}^*$

$$|H_m| \lesssim 1 \quad , \quad |\partial_x^k H_m| \lesssim 2^{k\kappa\sigma} \sum_{t'>0} 2^{(m-t')k} \mathcal{M}\mathcal{M}\Delta_{m-t'} f$$

(see Lemma 6.5). By the Leibniz rule, this implies

$$|\partial_x^{a+1}(\omega_{m-t} H_m)| \lesssim 2^{(m-t)(a+1)}\omega_{m-t} + 2^{(a+2)\kappa\sigma} \sum_{t'>0} \sum_{\ell=0}^{a} 2^{(t'-t)\ell} 2^{(a+1)(m-t')} \mathcal{M}\mathcal{M}\Delta_{m-t} f \cdot \mathcal{M}\mathcal{M}\Delta_{m-t'} f.$$

We split the sum over $t'$ into two parts according to whether $t' \leq t$ or $t' > t$ and we estimate by $\|f\|_{\dot{F}_q^{\alpha,p}}$ the factor $\mathcal{M}\mathcal{M}\Delta_{m-t} f$ or $\mathcal{M}\mathcal{M}\Delta_{m-t'} f$ respectively (here we use (2.12)):

$$\sum_{t'>0} \sum_{\ell=0}^{a} 2^{(t'-t)\ell} 2^{(a+1)(m-t')} \mathcal{M}\mathcal{M}\Delta_{m-t} f \cdot \mathcal{M}\mathcal{M}\Delta_{m-t'} f$$

$$\lesssim \|f\|_{\dot{F}_q^{\alpha,p}} \left( \sum_{0<t'\leq t} 2^{(a+1)(m-t')} \mathcal{M}\mathcal{M}\Delta_{m-t'} f + \sum_{t'>t} 2^{(t'-t)a} 2^{(a+1)(m-t')} \mathcal{M}\mathcal{M}\Delta_{m-t} f \right)$$

$$\lesssim \|f\|_{\dot{F}_q^{\alpha,p}} \sum_{0<t'\leq t} 2^{(a+1)(m-t')} \mathcal{M}\mathcal{M}\Delta_{m-t'} f.$$

Putting these together,

$$|\partial_x^{a+1}(\omega_{m-t} H_m)| \lesssim 2^{(m-t)(a+1)}\omega_{m-t} + 2^{(a+2)\kappa\sigma} \|f\|_{\dot{F}_q^{\alpha,p}} \sum_{0<t'\leq t} 2^{(a+1)(m-t')} \mathcal{M}\mathcal{M}\Delta_{m-t'} f.$$

Setting $B_{a,\alpha}(t) := \sum_{0<t'\leq t} 2^{(\alpha-1-a)t'}$, we deduce that

$$\|\|2^{(\alpha-1-a)m}(\partial_x^{a+1}(\omega_{m-t} H_m))\|_{\ell^q(m)}\|_{L^p}$$
$$\lesssim 2^{(\alpha-1-a)t} \|\|2^{\alpha(m-t)}\omega_{m-t}\|_{\ell^q(m)}\|_{L^p} + 2^{(a+2)\kappa\sigma} \|f\|_{\dot{F}_q^{\alpha,p}} \sum_{0<t'\leq t} 2^{-(a+1)t'} \|\|2^{\alpha m}\mathcal{M}\mathcal{M}\Delta_{m-t'} f\|_{\ell^q(m)}\|_{L^p}$$
$$\lesssim 2^{(\alpha-1-a)t} \|\|2^{\alpha m}\omega_m\|_{\ell^q(m)}\|_{L^p} + 2^{(a+2)\kappa\sigma} \|f\|_{\dot{F}_q^{\alpha,p}} B_{a,\alpha}(t) \|\|2^{\alpha m}\mathcal{M}\mathcal{M}\Delta_m f\|_{\ell^q(m)}\|_{L^p}$$

By Proposition 2.5, this implies

$$\|\|2^{(\alpha-1-a)m}(\partial_x^{a+1}(\omega_{m-t} H_m))\|_{\ell^q(m)}\|_{L^p}$$
$$\lesssim 2^{(\alpha-1-a)t} \|\|2^{\alpha m}\omega_m\|_{\ell^q(m)}\|_{L^p} + 2^{(a+2)\kappa\sigma} \|f\|_{\dot{F}_q^{\alpha,p}} B_{a,\alpha}(t) \|\|2^{\alpha m}\Delta_m f\|_{\ell^q(m)}\|_{L^p}$$
$$\leq 2^{\kappa\sigma} 2^{(\alpha-1-a)t} \|f\|_{\dot{F}_q^{\alpha,p}} + B_{a,\alpha}(t) 2^{(a+2)\kappa\sigma} \|f\|_{\dot{F}_q^{\alpha,p}}^2,$$

where the last line is a consequence of Proposition 4.6.

Coming back to (6.3), we get
(6.4)
$$\|\|2^{(\alpha-1)m}\Delta_m(\partial_x(\omega_{m+r-t} H_{m+r}))\|_{\ell^q(m)}\|_{L^p} \lesssim 2^{-(\alpha-1-a)r} \left( 2^{\kappa\sigma} 2^{(\alpha-1-a)t} \|f\|_{\dot{F}_q^{\alpha,p}} + B_{a,\alpha}(t) 2^{(a+2)\kappa\sigma} \|f\|_{\dot{F}_q^{\alpha,p}}^2 \right).$$

Choose now $a = [\alpha]$, so that $B_{a,\alpha}(t) \lesssim 1$. Summing up on $r \leq 0$, which is possible since $\alpha - a - 1 < 0$, we thus obtain

$$\sum_{r\leq 0} \|\|2^{(\alpha-1)m}\Delta_m(\partial_x(\omega_{m+r-t} H_{m+r}))\|_{\ell^q(m)}\|_{L^p} \lesssim 2^{\kappa\sigma} 2^{(\alpha-1-a)t} \|f\|_{\dot{F}_q^{\alpha,p}} + 2^{(a+2)\kappa\sigma} \|f\|_{\dot{F}_q^{\alpha,p}}^2.$$

Summing over $t$, we finally get

$$(6.5) \quad \sum_{\substack{t>0, \\ t\equiv 0(\mathrm{mod}\ R)}} 2^{-\alpha t} \sum_{r\leq 0} \|\|2^{(\alpha-1)m}\Delta_m(\partial_x(\omega_{r+m-t}H_{r+m}))\|_{\ell^q(m)}\|_{L^p}$$

$$\lesssim 2^{\kappa\sigma} 2^{-(a+1)R}\|f\|_{\dot{F}_q^{\alpha,p}} + 2^{(a+2)\kappa\sigma} 2^{-\alpha R}\|f\|_{\dot{F}_q^{\alpha,p}}^2.$$

### 6.3. Estimate of $\sum_{0\leq r\leq a_\alpha t}$.

Applying (6.4) with $a=0$, we obtain

$$\|\|2^{(\alpha-1)m}\Delta_m(\partial_x(\omega_{m+r-t}H_{m+r}))\|_{\ell^q(m)}\|_{L^p} \lesssim 2^{-(\alpha-1)r}\left(2^{\kappa\sigma} 2^{(\alpha-1)t}\|f\|_{\dot{F}_q^{\alpha,p}} + B_{0,\alpha}(t) 2^{2\kappa\sigma}\|f\|_{\dot{F}_q^{\alpha,p}}^2\right),$$

where

$$B_{0,\alpha}(t) = \begin{cases} C & \text{if } \alpha < 1, \\ Ct & \text{if } \alpha = 1, \\ C2^{(\alpha-1)t} & \text{if } \alpha > 1. \end{cases}$$

Summing on $0\leq r\leq a_\alpha t$, we get

$$\sum_{0\leq r\leq a_\alpha t} \|\|2^{(\alpha-1)m}\Delta_m(\partial_x(\omega_{m+r-t}H_{m+r}))\|_{\ell^q(m)}\|_{L^p} \lesssim A_\alpha(t)\left(2^{\kappa\sigma} 2^{(\alpha-1)t}\|f\|_{\dot{F}_q^{\alpha,p}} + B_{0,\alpha}(t) 2^{2\kappa\sigma}\|f\|_{\dot{F}_q^{\alpha,p}}^2\right),$$

where

$$A_\alpha(t) = \begin{cases} C2^{(1-\alpha)a_\alpha t} & \text{if } \alpha < 1, \\ Ca_\alpha t & \text{if } \alpha = 1, \\ C & \text{if } \alpha > 1. \end{cases}$$

Summing over $t$, one gets

$$\sum_{\substack{t>0 \\ t\equiv 0(\mathrm{mod}\ R)}} 2^{-\alpha t} \sum_{0\leq r\leq a_\alpha t} \|\|2^{(\alpha-1)m}\Delta_m(\partial_x(\omega_{m+r-t}H_{m+r}))\|_{\ell^q(m)}\|_{L^p}$$

$$\lesssim \sum_{t\geq R} 2^{-\alpha t} A_\alpha(t)\left(2^{\kappa\sigma} 2^{(\alpha-1)t}\|f\|_{\dot{F}_q^{\alpha,p}} + 2^{2\kappa\sigma} B_{0,\alpha}(t)\|f\|_{\dot{F}_q^{\alpha,p}}^2\right).$$

When $\alpha > 1$, one has $a_\alpha = 1$. Since $A_\alpha(t) = C$ and $B_{0,\alpha}(t) = C2^{(\alpha-1)t}$, this implies

$$\sum_{\substack{t>0, \\ t\equiv 0(\mathrm{mod}\ R)}} 2^{-\alpha t} \sum_{0\leq r\leq a_\alpha t} \|\|2^{(\alpha-1)m}\Delta_m(\partial_x(\omega_{m+r-t}H_{m+r}))\|_{\ell^q(m)}\|_{L^p} \lesssim 2^{\kappa\sigma} 2^{-R}\|f\|_{\dot{F}_q^{\alpha,p}} + 2^{2\kappa\sigma} 2^{-R}\|f\|_{\dot{F}_q^{\alpha,p}}^2.$$

When $\alpha = 1$, the choice $a_\alpha = 1$ leads to $A_\alpha(t) = Ct$ and $B_{0,\alpha}(t) = Ct$. Hence

$$\sum_{\substack{t>0, \\ t\equiv 0(\mathrm{mod}\ R)}} 2^{-\alpha t} \sum_{0\leq r\leq a_\alpha t} \|\|2^{(\alpha-1)m}\Delta_m(\partial_x(\omega_{m+r-t}H_{m+r}))\|_{\ell^q(m)}\|_{L^p} \lesssim 2^{\kappa\sigma} R2^{-R}\|f\|_{\dot{F}_q^{\alpha,p}} + 2^{2\kappa\sigma} R^2 2^{-R}\|f\|_{\dot{F}_q^{\alpha,p}}^2.$$

When $\alpha < 1$, $A_\alpha(t) = C2^{(1-\alpha)a_\alpha t}$ and $B_{0,\alpha}(t) = C$. One needs to take $a_\alpha < \frac{\alpha}{1-\alpha}$ for the sum to converge and one gets

$$\sum_{\substack{t>0,\\ t\equiv 0 \,(\text{mod } R)}} 2^{-\alpha t} \sum_{0\leq r\leq a_\alpha t} \|\|2^{(\alpha-1)m}\Delta_m(\partial_x(\omega_{m+r-t}H_{m+r}))\|_{\ell^q(m)}\|_{L^p}$$

$$\lesssim 2^{\kappa\sigma}2^{-R(1-(1-\alpha)a_\alpha)}\|f\|_{\dot{F}_q^{\alpha,p}} + 2^{2\kappa\sigma}2^{-R(\alpha-(1-\alpha)a_\alpha)}\|f\|^2_{\dot{F}_q^{\alpha,p}}.$$

In any case, for every $\alpha > 0$, and assuming further that $a_\alpha \leq \alpha$ when $\alpha < 1$, we have

(6.6) $$\sum_{\substack{t>0,\\ t\equiv 0 \,(\text{mod } R)}} 2^{-\alpha t} \sum_{0\leq r\leq a_\alpha t} \|\|2^{(\alpha-1)m}\Delta_m(\partial_x(\omega_{m+r-t}H_{m+r}))\|_{\ell^q(m)}\|_{L^p}$$

$$\lesssim 2^{\kappa\sigma}R2^{-R\min(1,\alpha a_\alpha)}\|f\|_{\dot{F}_q^{\alpha,p}} + 2^{2\kappa\sigma}R^2 2^{-R\min(1,\alpha a_\alpha)}\|f\|^2_{\dot{F}_q^{\alpha,p}}.$$

6.4. **Conclusion.** The three subsections above, namely inequalities (6.2), (6.5) and (6.6), imply the desired estimate (6.1). This completes the proof of Proposition 6.1.

7. COMPLETION OF THE PROOF OF PROPOSITION 3.3

Let us summarize the current state of the proof. For every $\sigma \in \mathbb{N}^*$, we define

$$R := \frac{\kappa + 1}{\min(1, \alpha a_\alpha)}\sigma$$

where $a_\alpha$ has been introduced in Proposition 6.1. This automatically implies the condition (4.19) since $\alpha R \geq (\kappa + 1)\sigma/a_\alpha > 1$. Then by Lemma 4.10, Proposition 5.1 and Proposition 6.1, there exists $\eta > 0$ such that for every $f \in \dot{F}_q^{\alpha,p}$ with $\|f\|_{\dot{F}_q^{\alpha,p}} \leq \eta 2^{-\kappa\sigma}$, there exists a map $F = \tilde{g} + \tilde{h}$ in $\dot{F}_q^{\alpha,p}$ such that

$$\|\tilde{g}\|_{L^\infty} \lesssim R, \qquad \|\tilde{h}\|_{L^\infty} \lesssim 1$$

and in the good directions $x'$:

$$\|\partial_{x'}(h-\tilde{h})\|_{\dot{F}_q^{\alpha-1,p}} \lesssim \sigma^3\left(2^{(-\min(1,\alpha)+\frac{\kappa}{p})\sigma}\|f\|_{\dot{F}_q^{\alpha,p}} + 2^{\left(\max(1-\alpha,0)+\kappa\left(1+[\alpha]+\frac{1}{p}\right)\right)\sigma}\|f\|^2_{\dot{F}_q^{\alpha,p}}\right),$$

while in all directions:

$$\|\partial_x(h-\tilde{h})\|_{\dot{F}_q^{\alpha-1,p}} \lesssim \sigma^3\left(2^{\frac{\kappa}{p}\sigma}\|f\|_{\dot{F}_q^{\alpha,p}} + 2^{\left(\max(1-\alpha,0)+\kappa\left(1+[\alpha]+\frac{1}{p}\right)\right)\sigma}\|f\|^2_{\dot{F}_q^{\alpha,p}}\right),$$

$$\|\partial_x(g-\tilde{g})\|_{\dot{F}_q^{\alpha-1,p}} \lesssim \left(\sigma 2^{-\sigma}\|f\|_{\dot{F}_q^{\alpha,p}} + \sigma^2 2^{([\alpha]+1)\kappa\sigma}\|f\|^2_{\dot{F}_q^{\alpha,p}}\right).$$

In order to prove Proposition 3.3, we take for every $\delta > 0$ an integer $\sigma > 0$ such that

$$\sigma^3 2^{(-\min(1,\alpha)+\frac{\kappa}{p})\sigma} \leq \frac{\delta}{2}, \quad \sigma 2^{-\sigma} \leq \frac{\delta}{2}.$$

This is possible in view of the fact that $\frac{\kappa}{p} < \min(1,\alpha)$. This implies that

$$\|\partial_{x'}(f-F)\|_{\dot{F}_q^{\alpha-1,p}} \leq \|\partial_{x'}(h-\tilde{h})\|_{\dot{F}_q^{\alpha-1,p}} + \|\partial_{x'}(g-\tilde{g})\|_{\dot{F}_q^{\alpha-1,p}} \leq \delta\|f\|_{\dot{F}_q^{\alpha,p}} + D_\delta\|f\|^2_{\dot{F}_q^{\alpha,p}}$$

for some $D_\delta > 0$ which may depend on $\sigma$ (and thus on $\delta$). We also have

$$\|F\|_{L^\infty} \leq \|\tilde{g}\|_{L^\infty} + \|\tilde{h}\|_{L^\infty} \leq D_\delta,$$

and using Proposition 2.3

$$\|F\|_{\dot{F}_q^{\alpha,p}} \leq \left\|h - \tilde{h}\right\|_{\dot{F}_q^{\alpha,p}} + \|g - \tilde{g}\|_{\dot{F}_q^{\alpha,p}} + \|f\|_{\dot{F}_q^{\alpha,p}}$$
$$\leq D_\delta,$$

by enlarging $D_\delta$ if necessary. This completes the proof of Proposition 3.3.

## 8. Solving Hodge systems

*Proof of Theorem 1.2:* We follow [5, p. 284]. On the space $\dot{F}_q^{s,p}(\Lambda^l \mathbb{R}^d)$ of $l$-forms with coefficients in $\dot{F}_q^{s,p}(\mathbb{R}^d)$, we use the norm

$$\|\lambda\|_{\dot{F}_q^{s,p}(\Lambda^l \mathbb{R}^d)} = \max_{|I|=l} \|\lambda_I\|_{\dot{F}_q^{s,p}(\mathbb{R}^d)}$$

if $\lambda = \sum_{|I|=l} \lambda_I dx_I$. Let $\varphi \in \dot{F}_q^{\alpha,p}(\Lambda^l \mathbb{R}^d)$. Since $d : \dot{F}_q^{\alpha,p}(\Lambda^l \mathbb{R}^d) \to \dot{F}_q^{\alpha-1,p}(\Lambda^l \mathbb{R}^d)$ is bounded with closed range, by the open mapping theorem, there exists $\lambda^{(0)} \in \dot{F}_q^{\alpha,p}(\Lambda^l \mathbb{R}^d)$ such that

$$d\lambda^{(0)} = d\varphi$$

and

(8.1) $$\left\|\lambda^{(0)}\right\|_{\dot{F}_q^{\alpha,p}(\Lambda^l \mathbb{R}^d)} \leq C \|d\varphi\|_{\dot{F}_q^{\alpha-1,p}(\Lambda^l \mathbb{R}^d)}.$$

Choose $\delta > 0$ such that $C\delta \leq \frac{1}{2}$. Let $I \subset \mathbb{N}^d$ be a multi-index with length $l$. Theorem 1.1 provides a function $\beta_I^{(0)} \in \dot{F}_q^{\alpha,p}(\mathbb{R}^d) \cap L^\infty(\mathbb{R}^d)$ such that, for all $j \in [\![1, d]\!] \setminus I$ (note that there are at most $\kappa$ such indexes),

$$\left\|\partial_j \left(\lambda_I^{(0)} - \beta_I^{(0)}\right)\right\|_{\dot{F}_q^{\alpha-1,p}(\mathbb{R}^d)} \leq \delta \left\|\lambda_I^{(0)}\right\|_{\dot{F}_q^{\alpha,p}(\mathbb{R}^d)} \leq C\delta \|d\varphi\|_{\dot{F}_q^{\alpha-1,p}(\Lambda^l \mathbb{R}^d)}$$

and

$$\left\|\beta_I^{(0)}\right\|_{\dot{F}_q^{\alpha,p}(\mathbb{R}^d)} + \left\|\beta_I^{(0)}\right\|_{L^\infty(\mathbb{R}^d)} \leq C_\delta \left\|\lambda_I^{(0)}\right\|_{\dot{F}_q^{\alpha,p}(\mathbb{R}^d)} \leq C'_\delta \|d\varphi\|_{\dot{F}_q^{\alpha-1,p}(\Lambda^l \mathbb{R}^d)}.$$

Set $\beta^{(0)} := \sum_I \beta_I^{(0)} dx^I$. Then, $\beta^{(0)} \in \dot{F}_q^{\alpha,p}(\Lambda^l \mathbb{R}^d) \cap L^\infty(\Lambda^l \mathbb{R}^d)$. Moreover, since $d\varphi = d\lambda^{(0)}$,

$$\left\|d(\varphi - \beta^{(0)})\right\|_{\dot{F}_q^{\alpha-1,p}(\Lambda^l \mathbb{R}^d)} = \max_{|I|=l} \max_{j \notin I} \left\|\partial_j (\lambda_I^{(0)} - \beta_I^{(0)})\right\|_{\dot{F}_q^{\alpha-1,p}(\mathbb{R}^d)} \leq \frac{1}{2} \|d\varphi\|_{\dot{F}_q^{\alpha-1,p}(\Lambda^l \mathbb{R}^d)}$$

and

$$\left\|\beta^{(0)}\right\|_{\dot{F}_q^{\alpha,p}(\Lambda^l \mathbb{R}^d)} + \left\|\beta^{(0)}\right\|_{L^\infty(\Lambda^l \mathbb{R}^d)} \leq C' \|d\varphi\|_{\dot{F}_q^{\alpha-1,p}(\Lambda^l \mathbb{R}^d)}.$$

The same argument, applied to $2(\varphi - \beta^{(0)})$ instead of $\varphi$, yields $\beta^{(1)} \in \dot{F}_q^{\alpha,p}(\Lambda^l \mathbb{R}^d) \cap L^\infty(\Lambda^l \mathbb{R}^d)$ such that

$$\left\|\beta^{(1)}\right\|_{\dot{F}_q^{\alpha,p}(\Lambda^l \mathbb{R}^d)} + \left\|\beta^{(1)}\right\|_{L^\infty(\Lambda^l \mathbb{R}^d)} \leq C' \left\|2d\left(\varphi - \beta^{(0)}\right)\right\|_{\dot{F}_q^{\alpha-1,p}(\Lambda^l \mathbb{R}^d)} \leq C' \|d\varphi\|_{\dot{F}_q^{\alpha-1,p}(\Lambda^l \mathbb{R}^d)}$$

and

$$\left\|d\varphi - d\left(\beta^{(0)} + \frac{1}{2}\beta^{(1)}\right)\right\|_{\dot{F}_q^{\alpha-1,p}(\Lambda^l \mathbb{R}^d)} \leq \frac{1}{4} \|d\varphi\|_{\dot{F}_q^{\alpha-1,p}(\Lambda^l \mathbb{R}^d)}.$$

Iterating this procedure, we construct a sequence $(\beta^{(i)})_{i\geq 0} \in \dot{F}_q^{\alpha,p}(\Lambda^l \mathbb{R}^d) \cap L^\infty(\Lambda^l \mathbb{R}^d)$ such that, for all $N \geq 0$,

$$\left\| d\varphi - d\left(\sum_{i=0}^N 2^{-i}\beta^{(i)}\right)\right\|_{\dot{F}_q^{\alpha-1,p}(\Lambda^l \mathbb{R}^d)} \leq \frac{1}{2^{N+1}} \|d\varphi\|_{\dot{F}_q^{\alpha-1,p}(\Lambda^l \mathbb{R}^d)}$$

and

$$\left\|\beta^{(N)}\right\|_{\dot{F}_q^{\alpha,p}(\Lambda^l \mathbb{R}^d)} + \left\|\beta^{(N)}\right\|_{L^\infty(\Lambda^l \mathbb{R}^d)} \leq C' \|d\varphi\|_{\dot{F}_q^{\alpha-1,p}(\Lambda^l \mathbb{R}^d)}.$$

Therefore, if $\psi := \sum_{i=0}^\infty 2^{-i}\beta^{(i)}$, $\psi$ satisfies all the conclusions of Theorem 1.2. $\square$

## 9. Appendix: some properties of Schwartz functions

As a consequence of Proposition 2.5,

**Lemma 9.1.** *There exists a constant $C$ which depends only on $\Delta$ such that for all $p, q \in (1, \infty)$, for all $(f_m)_{m \in \mathbb{Z}} \in L^p(\mathbb{R}^d; \ell^q(\mathbb{Z}))$, for all $k \in \mathbb{N}$,*

$$\|\|\partial_x^k \Delta_m f_m\|_{\ell^q(m)}\|_{L^p} \leq C \|\|2^{km} f_m\|_{\ell^q(m)}\|_{L^p}.$$

For every $\gamma = (\gamma_1, \ldots, \gamma_d) \in \mathbb{N}^d$, we denote

$$|\gamma| = \gamma_1 + \cdots + \gamma_d, \quad \partial^\gamma = \partial_1^{\gamma_1} \cdots \partial_d^{\gamma_d}.$$

Moreover, $X^\gamma$ is the polynomial $X_1^{\gamma_1} \ldots X_d^{\gamma_d}$ and for every polynomial $P(X_1, \ldots, X_d) = \sum_\gamma a_\gamma X^\gamma$, $P(D)$ is the differential operator $\sum_\gamma a_\gamma \partial^\gamma$.

**Lemma 9.2.** *Suppose $\phi$ is a Schwartz function on $\mathbb{R}^d$ and $m \in \mathbb{N}$. Assume that for every polynomial $P$ of degree less or equal to $m$, $\int_{\mathbb{R}^d} P(x)\phi(x)\,dx = 0$. Then for every $\gamma \in \mathbb{N}^d$ such that $|\gamma| = m+1$, there exists a Schwartz function $\phi^{(\gamma)}$ such that*

$$\phi = \sum_{|\gamma|=m+1} \partial^\gamma \phi^{(\gamma)}.$$

*Proof.* In terms of the Fourier transform $\psi$ of $\phi$, the assumption means that $\partial_\gamma \psi(0) = 0$ for every $|\gamma| \leq m$ while the conclusion amounts to the existence of Schwartz functions $\psi^{(\gamma)}$ such that

$$\psi = \sum_{|\gamma|=m+1} \xi^\gamma \psi^{(\gamma)}.$$

Let $\eta \in C_c^\infty(\mathbb{R}^d)$ is a smooth cut-off function with $\eta \equiv 1$ near the origin. By the Taylor formula:

$$\forall \xi \in \mathbb{R}^d, \quad \psi(\xi) = \sum_{|\gamma|=m+1} \xi^\gamma \left(\frac{1}{m!}\int_0^1 (1-t)^m \partial_\gamma \psi(t\xi)\,dt\right).$$

By the identity $|\xi|^{2(m+1)} = \sum_{|\gamma|=m+1} c_\gamma \xi^{2\gamma}$ with $c_\gamma = (m+1)!/(\gamma_1!\ldots\gamma_d!)$, we also have

$$\psi(\xi) = \eta(\xi) \sum_{|\gamma|=m+1} \xi^\gamma \left(\frac{1}{m!}\int_0^1 (1-t)^m \partial_\gamma \psi(t\xi)\,dt\right) + \frac{1-\eta(\xi)}{|\xi|^{2(m+1)}}\psi(\xi)\sum_{|\gamma|=m+1} c_\gamma \xi^{2\gamma}.$$

We can thus set

$$\psi^{(\gamma)}(\xi) := \eta(\xi)\left(\frac{1}{m!}\int_0^1 (1-t)^m \partial_\gamma \psi(t\xi)\,dt\right) + c_\gamma(1-\eta(\xi))\frac{\xi^\gamma}{|\xi|^{2(m+1)}}\psi(\xi).$$

The proof is complete.

□

## 10. Appendix: proof of Proposition 4.8

We begin with the following result.

**Proposition 10.1.** *Let $\{k_j\}_{j\in\mathbb{Z}}$ be a sequence of non-negative integrable functions on $\mathbb{R}^d$, with*

$$\sup_{j\in\mathbb{Z}} \|k_j\|_{L^1(\mathbb{R}^d)} \lesssim 1 \tag{10.1}$$

*and*

$$\int_{\|y\|\geq 4\|x\|} \sup_{j\in\mathbb{Z}} |k_j(y-x) - k_j(y)| dy \leq A \tag{10.2}$$

*for some constant $A \geq 1$. Then the associated maximal function*

$$\mathfrak{M}f := \sup_{j\in\mathbb{Z}} |f| * k_j \tag{10.3}$$

*is of weak-type $(1,1)$, and is bounded on $L^p(\mathbb{R}^d)$ for all $1 < p \leq \infty$; more precisely,*

$$\|\mathfrak{M}f\|_{L^{1,\infty}} \lesssim A\|f\|_{L^1} \tag{10.4}$$

*and*

$$\|\mathfrak{M}f\|_{L^p} \lesssim A^{1/p}\|f\|_{L^p}, \quad 1 < p \leq \infty. \tag{10.5}$$

*Also, $\mathfrak{M}$ satisfies a vector-valued weak-$L^1$ and strong-$L^p$ bound, namely*

$$\|\|\mathfrak{M}f_i\|_{\ell^q}\|_{L^{1,\infty}} \lesssim A\|\|f_i\|_{\ell^q}\|_{L^1}, \quad 1 < q \leq \infty, \tag{10.6}$$

*and*

$$\|\|\mathfrak{M}f_i\|_{\ell^q}\|_{L^p} \lesssim A^{1/p}\|\|f_i\|_{\ell^q}\|_{L^p}, \quad 1 < p \leq q \leq \infty. \tag{10.7}$$

The first part of the statement, namely (10.4) and (10.5), is essentially in the work of Zó [35], whose proof we reproduce below. One relies on a Banach-valued version of the singular integral theorem. To prove (10.4), we consider the Banach spaces $B_1 = \mathbb{C}$, $B_2 = \ell^\infty$, and the vector-valued singular integral

$$f \mapsto Tf := \{f * k_j\}_{j\in\mathbb{Z}}$$

which is a mapping of a $B_1$-valued function $f$ to a $B_2$-valued function $Tf = \{f * k_j\}_{j\in\mathbb{Z}}$. For technical reasons, we consider truncations of this operator $T$, namely $T_M f := \{f * k_j\}_{j\in\mathbb{Z}, |j|\leq M}$ for $M \in \mathbb{N}$, and show that the operator norm $\|T_M\|_{L^1 \to L^{1,\infty}(\ell^\infty)}$ is $\lesssim A$ uniformly for $M \in \mathbb{N}$. At almost every point $x \in \mathbb{R}^d$, the kernel $\{k_j(x)\}_{j\in\mathbb{Z}, |j|\leq M}$ can be thought of as a linear map from $B_1$ to $B_2$, whose operator norm is $\sup_{j\in\mathbb{Z}, |j|\leq M} |k_j(x)|$, and we note that the latter is in $L^1(\mathbb{R}^d)$ since

$$\int_{\mathbb{R}^d} \sup_{j\in\mathbb{Z}, |j|\leq M} |k_j(x)| dx \leq \int_{\mathbb{R}^d} \sum_{j\in\mathbb{Z}, |j|\leq M} |k_j(x)| dx \lesssim M. \tag{10.8}$$

Since $T_M \colon L^\infty \to L^\infty(\ell^\infty)$ with norm $\lesssim 1$, the vector-valued singular integral theorem ([2, Theorem 4.2]) gives

$$\|T_M\|_{L^1 \to L^{1,\infty}(\ell^\infty)} \lesssim A,$$

and interpolation in turn gives

$$\|T_M\|_{L^p \to L^p(\ell^\infty)} \lesssim A^{1/p}$$

for all $1 < p < \infty$. Both bounds being independent of $M$, we obtain (10.4) and (10.5) by letting $M \to +\infty$. Similarly, to prove (10.6), we consider the Banach spaces $B_1 = \ell^q$, $B_2 = \ell^q(\ell^\infty)$, and for each $M \in \mathbb{N}$ the vector-valued truncated singular integral

$$f = \{f_i\}_{i \in \mathbb{Z}} \mapsto T_M f := \{f_i * k_j\}_{i,j \in \mathbb{Z}, |j| \leq M}$$

which is a mapping of a $B_1$-valued function $f = \{f_i\}_{i \in \mathbb{Z}}$ to a $B_2$-valued function $T_M f = \{f_i * k_j\}_{i,j \in \mathbb{Z}, |j| \leq M}$. The $B_2$ norm of $T_M f$ at $x$ is by definition

$$\left( \sum_{i \in \mathbb{Z}} \sup_{\substack{j \in \mathbb{Z} \\ |j| \leq M}} |f_i * k_j(x)|^q \right)^{1/q}.$$

At almost every point $x \in \mathbb{R}^d$, the kernel $\{k_j(x)\}_{j \in \mathbb{Z}, |j| \leq M}$ can be thought of as a linear map from $B_1$ to $B_2$, whose operator norm is $\sup_{j \in \mathbb{Z}, |j| \leq M} |k_j(x)|$. As before, we note that this latter expression is in $L^1$ for all $M \in \mathbb{N}$. Now if $1 < q \leq \infty$ and $M \in \mathbb{N}$, (10.5) gives $T_M \colon L^q(\ell^q) \to L^q(\ell^q(\ell^\infty))$, with norm $\lesssim A^{1/q}$ uniformly in $M$. Hence the vector-valued singular integral theorem gives

$$\|T_M\|_{L^1(\ell^q) \to L^1(\ell^q(\ell^\infty))} \lesssim A,$$

and interpolation in turn gives

$$\|T_M\|_{L^p(\ell^q) \to L^p(\ell^q(\ell^\infty))} \lesssim A^{1/p}$$

for all $1 < p \leq q$. Therefore, (10.6) and (10.7) follow, as we let $M \to +\infty$. To apply Proposition 10.1, we use the following Lemma:

**Lemma 10.2.** *Suppose $\varphi \colon \mathbb{R}^d \to \mathbb{R}$ is a non-negative integrable function on $\mathbb{R}^d$ satisfying*

(10.9) $$\int_{\mathbb{R}^d} \varphi(y) dy \lesssim 1,$$

(10.10) $$\int_{\|y\| \geq R} \varphi(y) dy \lesssim R^{-1} \quad \text{for all } R \geq 1,$$

*and*

(10.11) $$\int_{\mathbb{R}^d} |\varphi(y-x) - \varphi(y)| dy \lesssim \|x\| \quad \text{for all } x \in \mathbb{R}^d \text{ with } \|x\| \leq 1.$$

*Define, for all $j \in \mathbb{Z}$ and all $x \in \mathbb{R}^d$,*

(10.12) $$\varphi_j(x) := 2^{jd} \varphi(2^j x),$$

*and define, for $r \in \mathbb{R}^d$,*

(10.13) $$k_j(x) := \varphi_j(x + 2^{-j} r).$$

*Then the kernels $k_j$ satisfy (10.2) with $A \lesssim \ln(2 + \|r\|)$, i.e.*

$$\int_{\|y\| \geq 4\|x\|} \sup_{j \in \mathbb{Z}} |k_j(y-x) - k_j(y)| dy \lesssim \ln(2 + \|r\|).$$

*Proof.* The proof is a variant of the argument in [23, Chapter II, Section 4.2]. Indeed, it suffices to replace the sup by a sum, and show that

$$\int_{\|y\|\geq 4\|x\|} \sum_{j\in\mathbb{Z}} |k_j(y-x) - k_j(y)| dy \lesssim \ln(2 + \|r\|). \tag{10.14}$$

We assume $\|r\| \geq 2$, for the case $\|r\| < 2$ follows from a simple modification of the following argument. We split the sum into three parts: $\sum_{j\in\mathbb{Z}} = \sum_{2^j\|x\|\leq 1} + \sum_{1<2^j\|x\|<\|r\|} + \sum_{2^j\|x\|\geq\|r\|}$. The first sum can be estimated using condition (10.11):

$$\int_{\|y\|\geq 4\|x\|} \sum_{2^j\|x\|\leq 1} |k_j(y-x) - k_j(y)| dy \leq \int_{\mathbb{R}^d} \sum_{2^j\|x\|\leq 1} 2^{jd} |\varphi(2^j(y-x)+r) - \varphi(2^j y + r)| dy$$

$$= \int_{\mathbb{R}^d} \sum_{2^j\|x\|\leq 1} |\varphi(y - 2^j x + r) - \varphi(y+r)| dy$$

$$\lesssim \sum_{2^j\|x\|\leq 1} 2^j \|x\|$$

$$\lesssim 1.$$

The second sum can be estimated using condition (10.9):

$$\int_{\|y\|\geq 4\|x\|} \sum_{1<2^j\|x\|<\|r\|} |k_j(y-x) - k_j(y)| dy \leq \sum_{1<2^j\|x\|<\|r\|} 2 \int_{\mathbb{R}^d} k_j(y) dy$$

$$= 2 \sum_{1<2^j\|x\|<\|r\|} \int_{\mathbb{R}^d} \varphi(y) dy$$

$$\lesssim \log\|r\|.$$

The last sum can be estimated using condition (10.10):

$$\int_{\|y\|\geq 4\|x\|} \sum_{2^j\|x\|\geq\|r\|} |k_j(y-x) - k_j(y)| dy \leq \sum_{2^j\|x\|\geq\|r\|} 2 \int_{\|y\|\geq 2\|x\|} k_j(y) dy$$

$$= \sum_{2^j\|x\|\geq\|r\|} 2 \int_{\|y\|\geq 2\cdot 2^j\|x\|} \varphi(y+r) dy$$

$$\lesssim \sum_{2^j\|x\|\geq\|r\|} \frac{1}{2^j\|x\|}$$

$$\lesssim \frac{1}{\|r\|} \lesssim 1.$$

Altogether we get (10.14). □

Combining Proposition 10.1 and Lemma 10.2, we see that:

**Corollary 10.3.** *If $\varphi$ is as in Lemma 10.2, $r \in \mathbb{R}^d$, and $\varphi_j$, $k_j$ and $\mathfrak{M}$ are as defined as in (10.12), (10.13) and (10.3) respectively, then $\mathfrak{M}$ is bounded on $L^p$ with norm $\lesssim [\ln(2 + \|r\|)]^{1/p}$, and it satisfies the vector-valued estimate*

$$\|\|\mathfrak{M}f_i\|_{\ell^q}\|_{L^p} \lesssim [\ln(2+\|r\|)]^{1/p} \|\|f_i\|_{\ell^q}\|_{L^p}, \quad 1 < p \leq q \leq \infty. \tag{10.15}$$

The above statement is reminiscent of similar estimates in the scalar-valued case [23, Chapter II.5.10] and [20, Theorem 4.1]. See also [14, Theorem 3.1] for a more general vector-valued estimate, that includes the case $p > q$.

We are now in position to prove Proposition 4.8. We apply Corollary 10.3 when $\varphi = T$, where $T(x) := \min(1, \|x\|^{-(d+1)})$ as in (3.2). We first verify conditions (10.9), (10.10) and (10.11) with $T$ in place of $\varphi$. Indeed (10.9) and (10.10) are obvious, and (10.11) follows since if $x \in \mathbb{R}^d$ with $\|x\| \leq 1$, then

$$\int_{\mathbb{R}^d} |T(y-x) - T(y)| dy \leq \int_{1-\|x\| \leq \|y\| \leq 1+\|x\|} 2 dy + \int_{\|y\| \geq 1+\|x\|} \left| \frac{1}{\|y-x\|^{d+1}} - \frac{1}{\|y\|^{d+1}} \right| dy$$

$$\lesssim \|x\| + \int_{\|y\| \geq 1+\|x\|} \frac{\|x\|}{\|y\|^{d+2}} dy$$

$$\lesssim \|x\|,$$

where the first inequality relies on the fact that $T(y-x) = T(y) = 1$ when $\|y\| \leq 1 - \|x\|$ while in the second inequality we have used the mean value theorem to estimate the integrand (if $\|y\| \geq 1 + \|x\|$ and $\|x\| \leq 1$, then $\|y\| \geq 2\|x\|$, so $\|y - tx\| \geq \|y\|/2$ for all $t \in [0,1]$, hence the desired estimate). Hence Corollary 10.3 applies. Now for any $r \in \mathbb{R}^d$ and any $j \in \mathbb{Z}$, we have $T_j|f_j|(x + 2^{-j}r) = |f_j| * k_j(x) \leq \mathfrak{M}f_j(x)$ for all $x \in \mathbb{R}^d$. Hence from Corollary 10.3, we conclude that

$$\|\|(T_j|f_j|(\cdot + 2^{-j}r))_{j \in \mathbb{Z}}\|_{\ell^q(\mathbb{Z})}\|_{L^p(\mathbb{R}^d)} \lesssim [\ln(2 + \|r\|)]^{\frac{1}{p}} \|\|(f_j)_{j \in \mathbb{Z}}\|_{\ell^q(\mathbb{Z})}\|_{L^p(\mathbb{R}^d)} \quad \text{if } 1 < p \leq q \leq \infty.$$

Duality between $L^p(\ell^q)$ and $L^{p'}(\ell^{q'})$ then shows that

$$\|\|(T_j|f_j|(\cdot + 2^{-j}r))_{j \in \mathbb{Z}}\|_{\ell^q(\mathbb{Z})}\|_{L^p(\mathbb{R}^d)} \lesssim [\ln(2 + \|r\|)]^{\frac{1}{p'}} \|\|(f_j)_{j \in \mathbb{Z}}\|_{\ell^q(\mathbb{Z})}\|_{L^p(\mathbb{R}^d)} \quad \text{if } 1 \leq q \leq p < \infty.$$

Since both $[\ln(2 + \|r\|)]^{\frac{1}{p}}$ and $[\ln(2 + \|r\|)]^{\frac{1}{p'}}$ are bounded from above by $\ln(2 + \|r\|)$, (4.14) then follows.

PIERRE BOUSQUET:
UNIVERSITÉ PAUL SABATIER, CNRS UMR 5219, INSTITUT DE MATHÉMATIQUES DE TOULOUSE
F-31062 TOULOUSE CEDEX 9, FRANCE
*E-mail address*: `pierre.bousquet@math.univ-toulouse.fr`

Emmanuel Russ:
Université Grenoble Alpes, CNRS UMR 5582, Institut Fourier
100 rue des Mathématiques, F-38610 Gières, France
*E-mail address*: emmanuel.russ@univ-grenoble-alpes.fr

Yi Wang:
John Hopkins University, Department of Mathematics
3400 N. Charles Street, 216 Krieger Hall, Baltimore, MD, 21218, USA
*E-mail address*: ywang@math.jhu.edu

Po Lam Yung:
The Chinese University of Hong-Kong, Department of Mathematics
Lady Shaw Building, Ma Liu Shui, Shatin, Hong Kong
*E-mail address*: plyung@math.cuhk.edu.hk